\documentclass[12pt]{article}
\usepackage[T1]{fontenc}
\usepackage{latexsym,amssymb}
\usepackage{amsmath}
\usepackage[mathscr]{eucal}
\usepackage{color}
\usepackage{xcolor}
\usepackage[abbrev]{amsrefs}

\newtheorem{Theorem}{\bf Theorem}[section]
\newtheorem{Lemma}{\bf Lemma}[section]
\newtheorem{Proposition}{\bf Proposition}[section]
\newtheorem{Corollary}{\bf Corollary}[section]
\newtheorem{Remark}{\bf Remark}[section]
\newtheorem{Example}{\bf Example}[section]
\newtheorem{Definition}{\bf Definition}[section]

\newenvironment{theorem}{\begin{Theorem}$\!\!\!$}{\end{Theorem}}
\newenvironment{lemma}{\begin{Lemma}$\!\!\!$}{\end{Lemma}}

\newenvironment{remark}{\begin{Remark}$\!\!\!$}{\end{Remark}}

\newenvironment{definition}{\begin{Definition}$\!\!\!$}{\end{Definition}}

\newcommand{\be}{\begin{equation}}
\newcommand{\ee}{\end{equation}}
\newcommand{\ba}{\begin{eqnarray}}
\newcommand{\ea}{\end{eqnarray}}
\newcommand{\bs}{\begin{equation*}}
\newcommand{\es}{\end{equation*}}

\newcommand{\Ginvext}{{\psi_f}}

\newcommand{\N}{\mathbb{N}}
\newcommand{\R}{\mathbb{R}}
\newcommand{\disp}{\displaystyle}

\newcommand{\dee}{{\rm{d}}}

\newcommand{\lam}{\lambda}
\newcommand{\pf}{{p_f}}
\newcommand{\pt}{{p_{\theta}}}
\newcommand{\qf}{{q_f}}

\newcommand{\CB}{{\mathcal{B}}}

\newcommand{\vp}{\varphi}

\numberwithin{equation}{section}

\def\Xint#1{\mathchoice
{\XXint\displaystyle\textstyle{#1}}%
{\XXint\textstyle\scriptstyle{#1}}%
{\XXint\scriptstyle\scriptscriptstyle{#1}}%
{\XXint\scriptscriptstyle\scriptscriptstyle{#1}}%
\!\int}
\def\XXint#1#2#3{{\setbox0=\hbox{$#1{#2#3}{\int}$}
\vcenter{\hbox{$#2#3$}}\kern-.5\wd0}}

\def\dashint{\Xint-}

\title{Local solvability and dilation-critical singularities of supercritical fractional heat equations}
\author{\qquad\\
Yohei Fujishima, Kotaro Hisa, Kazuhiro Ishige and Robert Laister}
\date{}
\begin{document}
\maketitle
\begin{abstract}
We consider
the Cauchy problem 
for fractional semilinear heat equations with supercritical nonlinearities and  
establish both necessary conditions and sufficient conditions for local-in-time solvability.
We introduce the notion of a {\it dilation-critical singularity} (DCS) of the initial data and show that such  singularities always exist for a large class of supercritical nonlinearities. 
Moreover, we provide exact formulae for such singularities.
\end{abstract}
\vspace{20pt}
\smallskip
\noindent Y.\,F.:  Department of Mathematical and Systems Engineering, Faculty of Engineering,\\
\qquad\qquad Shizuoka University, 3-5-1, Johoku, Hamamatsu, 432-8561, Japan.\\
\noindent
E-mail: {\tt fujishima@shizuoka.ac.jp}\\

\smallskip
\noindent K.\,H.:  Mathematical Institute, Tohoku University,\\
\qquad\qquad 6-3 Aoba, Aramaki, Aoba-ku, Sendai 980-8578, Japan.\\
\noindent
E-mail: {\tt kotaro.hisa.d5@tohoku.ac.jp}\\

\smallskip
\noindent K.\,I.:  Graduate School of Mathematical Sciences, The University of Tokyo,\\
\qquad\qquad 3-8-1, Komaba, Meguro-ku, Tokyo 153-8914, Japan.\\
\noindent
E-mail: {\tt ishige@ms.u-tokyo.ac.jp}\\

\smallskip
\noindent R.\,L.: School of Computer Science and Creative Technologies,\\
 \qquad\qquad University of the West of England, Bristol BS16 1QY, UK.\\
\noindent
E-mail: {\tt Robert.Laister@uwe.ac.uk}\\
\vspace{15pt}
\vspace{3pt}
\newpage
\section{Introduction}
\label{sec:Intro}
We consider the existence and nonexistence of local-in-time, nonnegative solutions to the Cauchy problem for the fractional semilinear  heat equation 
\begin{equation}
  \tag{SHE}
  \label{eq:SHE}
  \begin{cases}
    \partial_t u + (-\Delta)^\frac{\theta}{2}u =f(u), & x\in \mathbb{R}^N, \,\,\, t>0,
    \\[3pt]
    u(x,0) = \mu(x), & x\in \mathbb{R}^N,  
  \end{cases}
\end{equation}
where $\mu$ is a nonnegative and locally integrable function in ${\mathbb R}^N$.  
We assume throughout  that $N\ge 1$, $0<\theta\le 2$ and 
$f\colon [0,\infty)\to [0,\infty)$  is  continuous. 
The local solvability of problem~\eqref{eq:SHE} depends crucially upon the growth of $f(u)$ as $u\to\infty$ and 
the strength of any singularities in the initial data $\mu$. 
Typical nonlinearities  one may keep in mind, are
\begin{equation}
  \label{eq:power}
  f(u)=u^p,\qquad p>1,
\end{equation}
and
\begin{equation}
  \label{eq:1.1}
  f(u)=\exp(u^p),\qquad   p>0.
\end{equation}
In particular we will be concerned with {\it supercritical} nonlinearities $f$ (see Definition~\ref{def:pf}). 
Roughly speaking, one may think of a supercritical nonlinearity as one having a larger growth rate (in some sense) than the function $u\mapsto u^{\pt}$ at infinity,
where $\pt$ is a critical exponent given by 
\begin{equation}
\label{eq:critexp}
\pt:=1+\frac{\theta}{N}.
\end{equation}

Our primary  goal is to determine necessary conditions and sufficient conditions  for the  local  solvability of problem~\eqref{eq:SHE} 
for a  large class of monotone nonlinearities $f$. As a consequence, we are  able to  identify 
the critical strength of any isolated singularities of the initial data $\mu$ for  local solvability.
Results of this type are important for the development of any general function-analytic framework for well-posedness of problem~\eqref{eq:SHE}, 
for example in identifying an optimal Orlicz class of initial data for a given nonlinearity $f$. 
Such studies exist for a variety of semilinear problems: see  
\cites{FHIL, FI1} for problem~\eqref{eq:SHE} with power law and power-log nonlinearities; 
see \cites{FI01,FI02} for semilinear parabolic systems; 
\cite{HI02} for the linear heat equation with  nonlinear boundary conditions;
\cites{HS, HT} for the Hardy parabolic equation;
 \cite{IKO} for higher-order semilinear parabolic equations; 
 \cite{HIT} for semilinear heat equations in a half-space of ${\mathbb R}^N$. 

\subsection{Background} 
\label{subsec:background}
In all that  follows, we denote by $L^{1,+}_{{\rm loc}}(\R^N)$ the set of nonnegative and locally integrable functions in ${\mathbb R}^N$. 
The open Euclidean ball in $\mathbb{R}^N$ centered at $x$ with radius $r$, is denoted by $B(x,r)$; when $x=0$ we simply write $B_r$.  
For any  $E\subset{\mathbb R}^N$, we denote by $\displaystyle{\chi_E}$ the characteristic function of $E$; when $E=B_r$ we simply write $\chi_r$.

%
Consider first the special case of a 
power law nonlinearity $f(u)=u^p$, with $\mu\in L^{1,+}_{{\rm loc}}(\R^n)$:
\begin{equation}
 \tag{P}
 \label{eq:P}
  \begin{cases}
    \partial_t u + (-\Delta)^\frac{\theta}{2}u =u^p, & x\in \mathbb{R}^N, \,\,\, t>0,
    \\[3pt]
    u(x,0) = \mu(x), & x\in \mathbb{R}^N.
  \end{cases}
\end{equation}
The question of local  solvability of problem~\eqref{eq:P}   has been studied in several papers;  
see  \cites{BC, BP, FI, GM, HI01, KY, RS, W1, W2} and the references therein. 
See also \cite{QS}, which includes a comprehensive survey of  the literature pertaining to problem~\eqref{eq:P} when $\theta =2$.
Among others, in \cite{HI01} the  authors  studied qualitative properties of the initial trace
of nonnegative solutions to problem~\eqref{eq:P} and obtained necessary conditions and sufficient conditions for the existence of solutions. 
In particular, they proved the following. 
\begin{itemize}
  \item[(P1)] (Subcritical)
   Let $1<p<{\pt}$. 
  Problem~\eqref{eq:P} is locally solvable 
  if and only if 
  \[
  \sup_{x\in{\mathbb R}^N}\int_{B(x,1)} \mu(y)\, \dee y <\infty .
  \]
  \item[(P2)] (Critical)
   Let $p= {\pt}$  and 
  set $\mu_0(x)=\disp{|x|^{-N}\displaystyle{\biggl[\log\left(e+|x|^{-1} \right)\biggr]^{-\frac{N}{\theta}-1}}}$. 
  \begin{itemize}
    \item[(a)] There exists  $\gamma=\gamma  (N,\theta)>0$ such that if $\mu(x)\ge\gamma\mu_0(x)$
    for almost all (a.a.) $x$ in a neighborhood of the origin, then problem~\eqref{eq:P} is not locally solvable.
    \item[(b)]  There exists $\varepsilon =\varepsilon (N,\theta) >0$ such that if $\mu(x)\le \varepsilon \mu_0(x)+K$ 
    for a.a.~$x\in {\mathbb R}^N$, for some $K>0$, then problem~\eqref{eq:P} is locally solvable.
  \end{itemize}
  \item[(P3)] (Supercritical)  Let $p> {\pt}$ and 
   set  $\mu_0(x)=\displaystyle{|x|^{-\frac{\theta}{p-1}}}$. 
  \begin{itemize}
    \item[(a)] There exists  $\gamma=\gamma  (N,\theta ,p)>0$ such that if $\mu(x)\ge\gamma\mu_0(x) $ 
    for  a.a.~$x$ in a neighborhood of the origin, then problem~\eqref{eq:P} is not locally solvable.
    \item[(b)]  There exists $\varepsilon =\varepsilon (N,\theta ,p) >0$ such that if $\mu(x)\le \varepsilon \mu_0(x)+K$ 
    for a.a.~$x\in {\mathbb R}^N$, for some $K>0$, then problem~\eqref{eq:P} is locally solvable.
  \end{itemize}
\end{itemize}
The results in (P3)  reveal two things in the supercritical case. First, that the strength of the singularity  at the origin of the function
\[
\mu_0(x)=
 |x|^{-\frac{\theta}{p-1}}, \qquad x\in\R^N, 
\]
is in some sense a threshold for the local solvability of problem~\eqref{eq:P}. 
For example, if the initial data $\mu$ satisfies $\mu_0 =o(\mu)$ uniformly near the origin, 
then by (P3a) problem~\eqref{eq:P} is not locally solvable; 
conversely, if $\mu =o(\mu_0)$ uniformly near the origin and $\mu$ is bounded away from the origin, 
then by (P3b) problem~\eqref{eq:P} is locally solvable.   
Secondly,  within the parameterised class of initial data $\mu =\lam \mu_0$ ($\lam >0$) both existence and nonexistence can occur: 
existence  holds for $\lam$ small enough ($\lam\le\varepsilon $), while nonexistence pertains for $\lam$ large enough ($\lam\ge\gamma $). 
In fact one can easily show via supersolution methods that there exists a unique $\lam_0>0$ such that, 
with  $\mu =\lam \mu_0$, problem~\eqref{eq:P} is locally solvable for all $\lam\in (0,\lam_0 )$ while not  locally solvable for any $\lam > \lam_0$. 
Alternatively, setting $\mu_c=\lam_0\mu_0$  
we see that problem~\eqref{eq:P} with $\mu =\lam \mu_c$ is locally solvable for all $\lam\in (0,1 )$ and not locally solvable for any $\lam > 1$.
This type  of parameterisation, along the ray $\{ \lam \mu_c\}_{\lam >0}$, 
was considered in \cites{HIT} for problem~\eqref{eq:P} on a half-space, 
where $\mu_c$ is referred to as an `optimal singularity'. 
See also \cites{FHIL,HT}.

Equivalently, if one considers instead the parameterised family of {\it dilations} of $\mu_c$, namely  $\{  \mu_c(x/\lam)\}_{\lam >0}$, then we see from the above that problem~\eqref{eq:P} with $\mu (x) =\mu_c(x/\lam)$ is locally solvable for all $\lam\in (0,1 )$ and not  locally solvable for any $\lam > 1$.  
Accordingly (see Definition~\ref{def:CS}), we will term such a function $\mu_c$ a {\it dilation-critical singularity} (DCS). One of the key findings of this paper is that this  parameterised family of dilations is in fact the {\it natural} choice, even for the general semilinear problem~\eqref{eq:SHE} where homogeneous scaling properties of $f$ (such as exist for problem~\eqref{eq:P}) are not available. See Theorem~\ref{thm:main} and the preamble to Definition~\ref{def:CS}.  
%

Our first main result  (Theorem~\ref{thm:main}) generalises the existence and nonexistence results in (P3) to a general class of monotone, convex, supercritical nonlinearities $f$ for problem~\eqref{eq:SHE}. 
Our second main result  (Theorem~\ref{thm:CS}) is that for  $f$ in this class, 
problem~\eqref{eq:SHE} always possesses   dilation-critical singularities, for which we also provide exact  formulae. 

\medskip
In several works 
(see e.g.,  \cites{FI, GM, IJMS, Io, IoT, MS, RT, V}),  the case of an exponential-power nonlinearity $f(u)=\exp (u^p)$ ($p>0$) has been considered: 
\begin{equation}
  \tag{E}
  \label{eq:E}
  \begin{cases}
    \partial_t u + (-\Delta)^\frac{\theta}{2}u =\exp (u^p), & x\in \mathbb{R}^N, \,\,\, t>0,
    \\[3pt]
    u(x,0) = \mu(x)\ge 0, & x\in \mathbb{R}^N.
  \end{cases}
\end{equation}
Among others, in \cite{FI} 
problem~\eqref{eq:E} was studied  with $\theta=2$ and  $p\ge 1$ 
and  the following results were obtained (see \cite{FI}*{Theorems~1.1, 1.2 and 5.2}): 
\begin{itemize}
  \item[(E1)] 
  problem~\eqref{eq:E} is locally solvable
  if 
  $\disp{\left(\mu^{p-1}\exp(\mu^p)\right)^r\in L^{1}_{{\rm uloc}}\left({\mathbb R}^N\right)}$ with $ r>N/2$; 
  \item[(E2)] 
  problem~\eqref{eq:E} is locally solvable
   if 
  $\disp{\left(\mu^{p-1}\exp(\mu^p)\right)^{r}\in {\mathcal L}^{1}_{{\rm uloc}}\left({\mathbb R}^N\right)}$
  with   $r=N/2$;
  \item[(E3)]  for $0<r<N/2$, there exists $\mu\in L^{1,+}_{{\rm loc}}\left(\R^N\right)$ such that 
  $\disp{\left(\mu^{p-1}\exp(\mu^p)\right)^r\in L^{1}_{{\rm uloc}}\left({\mathbb R}^N\right)}$,
  but for which 
  problem~\eqref{eq:E} is not locally solvable.
\end{itemize}
Here, $L^1_{{\rm uloc}}\left({\mathbb R}^N\right)$ denotes the space of uniformly locally integrable functions  
and ${\mathcal L}^1_{{\rm uloc}}\left({\mathbb R}^N\right)$ 
is the closure  of the space of bounded uniformly continuous functions in ${\mathbb R}^N$  
with respect to the norm in ${L^1_{{\rm uloc}}\left({\mathbb R}^N\right)}$.
See also \cite{V} for sign-changing solutions when $p=1$ and $\theta =2$. 
For the fractional case $0<\theta <2$, see \cite{GM}*{Theorem~1.6} for (E1) and (E3)
when $p\ge 1$.  

It transpires that assertions~(E1)--(E3) are not enough to identify 
the DCS 
 of the initial data for the  local solvability of problem~(E) (see Remark~\ref{Remark:4.1}). 
However, as an application of our general results for the semilinear problem~\eqref{eq:SHE},  
we obtain refined results on the  local  solvability of problem~\eqref{eq:E}, 
together with a generalisation to a finite number of compositions of the exponential function  (Theorem~\ref{Theorem:GE}). 

\vspace{5pt}

\subsection{Main results}
Motivated by \cites{FI, GM, HI01}, 
we introduce the following assumptions on the nonlinearity:
\begin{itemize}
  \item[{(M)}] (Monotonicity) $f\colon [0,\infty)\to [0,\infty)$ is continuous, 
  nondecreasing and $f(u )>0$ for all $u >0$;
  \item[(S)] (Superlinearity) there exists $  \tau_0>0$ such that $f\in C^1([  \tau_0,\infty))$ and
  ${\displaystyle{ \int_{  \tau_0}^\infty \frac{\dee s }{f(s)} < \infty}}$; 
  \item[(C)] (Convexity) there exists $\tau_1\ge 0$ such that  $f$ is convex in $[  \tau_1,\infty)$.
\end{itemize}
For $f$ satisfying (M) and (S), we define 
\begin{equation}
  \label{eq:1.3}
  F(u)=\int_u^\infty \frac{\dee s}{f(s)} \quad\text{for}\quad u >0\qquad\text{and}\qquad F_0=\lim_{u\to 0}F(u),
\end{equation}
noting that $F_0$ may be finite or infinite. Now define
\begin{equation}
\label{eq:G}
G(u)=
\left\{
\begin{array}{ll}
\displaystyle{\frac{1}{F(u)}}, &  u >0,
\vspace{7pt}\\
\displaystyle{\frac{1}{F_0}}=:G_0, & u =0,
\end{array}
\right.
\end{equation}
with the convention that $G_0=0$ if $F_0=\infty$. 

The following hypothesis  plays a central role in this paper, necessary for  our notion of supercriticality  in Definition~\ref{def:pf}.

\begin{itemize}
\item[(L)]  (Limit) $\displaystyle{q_f:=\lim_{u\to\infty}f'(u)F(u)}$ exists and is finite. 
\end{itemize}

\begin{remark}\label{rem:ORV}
{\rm   
Suppose $f$ satisfies {\rm (M)} and {\rm (S)}.
\begin{itemize}
\item[{\rm (i)}] 
 If {\rm (L)} holds then  it follows from \cite{FI}*{Remark 1.1} that ${\qf}\ge 1$. Furthermore, 
\begin{equation} \notag \frac{1}{2q_f}f'(u) \le G(u)=\frac{1}{F(u)} \le \frac{2}{q_f} f'(u) \end{equation}
for large enough $u>0$. 
For sufficiently smooth $f$ (eventually $C^2$) it also follows from l'H\^opital's rule that
\[
  {\qf}=\lim_{u\to\infty}\frac{F(u)}{1/f'(u)}=\lim_{u\to\infty} \frac{f'(u)^2}{f(u)f''(u)}, 
\]
whenever the latter limit exists.
\item[{\rm (ii)}]
 If $f(0)>0$ then $G_0>0$. If $f(0)=0$ and $f$ is Lipschitz at zero, then $G_0=0$.
\item[{\rm (iii)}]   $F\colon (0,\infty )\to (0,F_0 )\ $  is strictly decreasing, with inverse $F^{-1}\colon (0,F_0 )\to (0,\infty )$.
\item[{\rm (iv)}]  $G\colon [0,\infty )\to [G_0,\infty )$  is strictly increasing, with inverse $G^{-1}\colon [G_0,\infty )\to [0,\infty )$. 
It will be useful for us to extend the domain of $G^{-1}$ 
via
$\Ginvext \colon [0,\infty )\to [0,\infty )$
as 
\begin{equation}
\label{eq:Ginv}
  \Ginvext(u) 
  =
  \begin{cases}
    0, & 0\le u\le G_0,\\
    G^{-1}(u), &  u> G_0.
  \end{cases}
\end{equation}
In particular, if  $G_0=0$ then $\Ginvext=G^{-1}$.
\end{itemize}
}
\end{remark}
\vspace{3mm}

It is via the function~$G$ that we study the critical threshold for the local solvability of problem~\eqref{eq:SHE} 
and subsequently obtain the dilation-critical singularity $\mu_c$  for  local  solvability. 

Following \cite{DF1}, and noting Remark~\ref{rem:ORV} (i), we now introduce an exponent measuring the growth of $f$ at infinity. As such we are able  to define precisely the meaning of the term `supercritical nonlinearity'. 

\begin{definition} 
\label{def:pf}
Suppose $f$ satisfies {\rm (M)}, {\rm (S)} and {\rm (L)}. 
We define the {\rm growth exponent $\pf\in (1,\infty]$ of $f$} to be the H\"{o}lder conjugate of $\qf$; i.e.,  $\pf$ satisfies
\begin{equation}
  \frac{1}{\pf}+\frac{1}{\qf}=1.
  \label{eq:pf}
\end{equation}
We then say that $f$ is {\rm\bf supercritical} if $\pf >\pt$ $($recalling \eqref{eq:critexp}$)$.
\end{definition}

We will abuse notation slightly and continue to use  $\pt $ for the critical exponent; i.e.,  where $\theta$ is a real number, not a function.

\begin{remark}
  \label{Remark:1.1}
  {\rm 
  We list some typical examples of nonlinearities satisfying (M), (S), (C) and (L), together with their growth rate exponents.
  \begin{itemize}
    \item[{(i)}] $f(u)=u^p$, $p>1$; $\pf =p$; supercritical if $p>\pt$.
    \item[{(ii)}] $f(u)=u^p+u^q$, $p>q>1$; $\pf =p$; supercritical if $p>\pt$.
    \item[{(iii)}] $f(u)=u^p\left(\log (1+u)\right)^q$, $p>1$, 
    $q>0$; $\pf =p$; supercritical if $p>\pt$.
    \item[{(iv)}] $f(u)=\exp_n(u^p)$, $p>0$, $n\in\N$, where $\exp_n$ denotes $n$ compositions of $\exp$; $\pf =\infty$.
  \end{itemize}
  }
\end{remark}
\vspace{3mm}

We introduce the  set $\CB$ of functions which are bounded away from the origin: 
\[
  \CB := \bigcap_{r>0} L^{\infty }\left(\R^n\backslash B_r\right). 
\]
Recalling \eqref{eq:G}, we are now in a position to state our first main result.


%
\begin{Theorem}[Local Solvability]
  \label{thm:main}
  Let  $\mu\in L^{1,+}_{{\rm loc}}(\R^N)$ and suppose $f$ satisfies {\rm (M)}, {\rm (S)} and {\rm (L)}. 
  \begin{itemize}
    \item[{\rm (a)}] If $f$ satisfies  {\rm (C)},  
    then there exists $\gamma=\gamma(N,\theta,f)>0$ 
    such that if 
    \begin{equation}
      \label{eq:1.5a}
      G(\mu(x))
      \ge
      \gamma |x|^{-\theta}
     \end{equation}
    for a.a.~$x$ in a neighborhood of the origin, 
    then problem~\eqref{eq:SHE} is not  locally solvable.
    \item[{\rm (b)}]  Let $\mu\in\CB$. If  $f$ is supercritical, 
    then there exists $\varepsilon=\varepsilon(N,\theta,f)>0$  such that if 
    \begin{equation}
      \label{eq:1.5}
      G(\mu(x))
      \le
      \varepsilon |x|^{-\theta} 
    \end{equation}
     for a.a.~$x$ in a neighborhood of the origin,
    then problem~\eqref{eq:SHE} is locally solvable.
    In particular, there exists a solution $u\in  L_{\rm{loc}}^\infty\left((0,T),L^\infty \left(\R^N\right)\right)$ for some $T>0$.
  \end{itemize}
\end{Theorem}

By `locally solvable'  we mean that problem~\eqref{eq:SHE}
possesses a nonnegative, local-in-time solution  (see Definition~\ref{Definition:2.1}). 
Under additional regularity assumptions on $f$, such as local Lipschitz continuity, 
the bounded solution guaranteed to exist in Theorem~\ref{thm:main}~(b), 
will be a local-in-time classical solution. 
As such one might wish to consider whether this solution can be continued globally in time. 
The behaviour of $f$ near zero then becomes the determining factor; see for example \cites{FI3,LS}.

To prove  Theorem~\ref{thm:main}~(a), we first obtain an {\it a priori} bound for any local solution of problem~\eqref{eq:SHE}, under minimal regularity assumptions on the solution (see Theorem~\ref{thm:nec}). The proof of Theorem~\ref{thm:nec}  follows that of \cite{FI}, itself a generalisation of the classical estimates in \cites{Fuj,W1}.

To prove Theorem~\ref{thm:main}~(b), we first obtain sufficient conditions for the  existence of a local supersolution  of problem~\eqref{eq:SHE} (see Theorem~\ref{thm:suff}). 
This involves a more delicate construction of a supersolution to problem~\eqref{eq:SHE} than those in \cites{FI, HI01, RS}. 
We adapt the proofs of Assertions~(E1) and (E2) in \cites{FI, GM} (based on  \cite{W1}),
where the following property holds:
\begin{itemize}
  \item[$\bullet$] 
  for large enough $u>0$,  $\displaystyle{ G(u)^{\frac{N}{\theta}}}$ is convex
  and $\displaystyle{G(\mu)^{\frac{N}{\theta}}}$ is locally integrable in ${\mathbb R}^N$.
\end{itemize}
Here however,  condition~\eqref{eq:1.5} is not enough to guarantee that 
$\displaystyle{G(\mu)^{\frac{N}{\theta}}}$ is locally integrable in ${\mathbb R}^N$
and  we cannot apply the arguments in \cites{FI, GM}. 
Instead we obtain a supersolution
under  the following property:
\begin{itemize}
  \item[$\bullet$] 
 for some $\beta >0$ and for all large enough $u>0$, $\displaystyle{G(u)^{\beta}}$
  is convex   and $\displaystyle{G(\mu)^{\beta}}$ is locally integrable in ${\mathbb R}^N$.
\end{itemize}

\medskip
It is clear from Theorem~\ref{thm:main} that the  local  solvability of problem~\eqref{eq:SHE}
depends crucially on the relative ordering,  as $x\to 0$, of $G(\mu)$ and the parameterised family of functions
\[
x\mapsto k |x|^{-\theta}=\left|\frac{x}{k^{\frac{1}{\theta}}}\right|^{-\theta},\qquad  k >0.
\] 
In order to characterise this behaviour more precisely,  we  introduce the following definition (replacing $k^{1/\theta}$ by $\lam$).
\begin{Definition} 
  \label{def:CS}
  We say that  $\mu_c\in \CB\cap L^{1,+}_{{\rm loc}}\left({\mathbb R}^N\right) $   
  is a {\rm{\bf dilation-critical singularity}  (DCS)}
  for  problem~\eqref{eq:SHE} if, 
  with $\mu(x) = \mu_c(x/\lambda)$, 
  the following hold: 
  \begin{itemize} 
    \item[{\rm (i)}]   for all $\lambda\in (0,1)$ 
    problem~\eqref{eq:SHE} is locally solvable; 
    \item[{\rm (ii)}] for all $\lambda >1$ 
    problem~\eqref{eq:SHE} is  not locally solvable. 
  \end{itemize} 
\end{Definition}

We  recall  some asymptotic notation. For any nonnegative functions $\phi$ and  $\psi$ 
we write $\phi\asymp \psi$  whenever there exist constants $C_1, C_2>0$ and $L>0$ such that $C_1\phi (u)\le \psi (u)\le C_2\phi (u)$ for all $u\ge L$. We may now state our second main result.


\begin{Theorem}[Dilation-Critical Singularity]
  \label{thm:CS}
  Suppose $f$ satisfies {\rm (M)}, {\rm (S)}, {\rm (C)}, {\rm (L)} and is supercritical. 
 Let $\vp ,\psi\colon [0,\infty )\to [0,\infty ) $
 be any nondecreasing functions satisfying $\vp\asymp G$ {\rm (}or equivalently, by $($L$)$,  $\vp\asymp f'${\rm )} 
 and $\vp\circ\psi\asymp\mathrm{id}$. 
 Then there exists a $unique$ $\lambda_0=\lambda_0 (\psi) >0$ such that 
  \begin{equation}
    \label{eq:critsing}
    \mu_c(x):= \psi\left(
     \lambda_0|x|^{-\theta}\right),\qquad x\in\R^N, 
  \end{equation}
  is a dilation-critical singularity 
  for problem~\eqref{eq:SHE}.
\end{Theorem}
Taking $\vp =G$ and $\psi =\Ginvext$ in Theorem~\ref{thm:CS},  immediately yields the following:

\begin{Corollary}
\label{cor:CS} 
Suppose $f$ satisfies {\rm (M)}, {\rm (S)}, {\rm (C)}, {\rm (L)} and is supercritical. 
Then there exists a unique ${\lambda}_0 >0$ such that 
\begin{equation}
\label{eq:GCS}
\mu_c(x):= 
\Ginvext 
\left(
{\lambda}_0|x|^{-\theta}\right),\qquad x\in\R^N, 
\end{equation}
is a dilation-critical singularity 
for problem~\eqref{eq:SHE}, with $\Ginvext$ as defined in \eqref{eq:Ginv}.
\end{Corollary}

\begin{Remark}\label{rem:asym}
{\rm 
As per Corollary~\ref{cor:CS},  functions $\vp$ and $\psi$ satisfying the hypotheses of Theorem~\ref{thm:CS}
always exist.
However, in applications it may be difficult  to compute explicitly  the function  $\Ginvext$; 
indeed, the power of Theorem~\ref{thm:CS}  lies in the potential to apply asymptotic methods to reveal the nature of a DCS. 
See Example~\ref{eg:powerlog}  and Theorem~\ref{Theorem:GE}(c).
}
\end{Remark}

We briefly illustrate some applications of Theorem~\ref{thm:CS} in some simple or known cases.  
In the first two examples $G^{-1}$ is explicitly computable so that we may use \eqref{eq:GCS}; 
in the third example, the DCS is obtained via asymptotic means utilising \eqref{eq:critsing}.


\begin{Example}
\label{eg:power}
{\rm 
Let $f(u)=u^p$, where $p>\pt$. 
By Remark~\ref{rem:ORV}~(ii) we have $G_0=0$ and $\Ginvext (u)=G^{-1}(u)=c_pu^{1/(p-1)}$, 
where $c_p=(p-1)^{-1/(p-1)}$. By Corollary~\ref{cor:CS} 
a DCS is  given by
 \[
 \mu_c(x)=\hat{\lambda}_0 |x|^{-\theta/(p-1)}
 ,\qquad x\in\R^N,
 \]
 for some unique  $\hat{\lam}_0>0$. See \cite{HI01} and the discussion around (P3) in Section~\ref{subsec:background} above.
}
\end{Example}

\begin{Example}
\label{eg:exp}
{\rm 
Let $f(u)=\exp (u)$. We have $G_0=1$ and $G^{-1}(u)=\log u$, for $u\in [1,\infty)$. 
Hence, by Corollary~\ref{cor:CS}, there exists  ${\lam}_0 >0$ such that 
 \[
 \mu_c(x)=
 \Ginvext
 (\lambda_0 |x|^{-\theta})
 =
 \log (  {\lambda}_0 |x|^{-\theta}) \chi_r(x),
 \qquad x\in\R^N,
 \]
is a DCS, where $r={\lambda}_0^{1/\theta}$.
}
\end{Example}

\begin{Example}
\label{eg:powerlog}
{\rm 
Let 
$f_0(u)=u^p\left(\log u\right)^q$ for $u\ge L>1$,
where $p>\pt$ and $q\in{\mathbb R}$. For $L$  large enough, $f_0$ is strictly increasing and convex. Moreover, 
\[ f_0'(u)\asymp u^{p-1}\left(\log u \right)^q=:\varphi_0 (u)\]
and  $\vp_0$ is strictly increasing on $[L,\infty )$ for $L$ large enough. Furthermore, by \cite{FHIL}*{Lemma 2.6}, 
\begin{equation}\label{eq:ORV1}
 \vp_0^{-1}(u)\asymp u^{\frac{1}{p-1}}\left(\log u\right)^{-\frac{q}{p-1}}=:\psi_0 (u),
\end{equation}
and  $\psi_0$ is strictly increasing on $[L,\infty )$ for $L$ large enough.
 Let us denote the extensions by zero of  $\vp_0$ and $\psi_0$  to the domain $[0,\infty )$ by $\vp$ and $\psi$, respectively.
Now take $f\colon [0,\infty )\to [0,\infty )$ to be any  nondecreasing continuous function  satisfying $f>0$ in $(0,L)$ and $f=f_0$ in $[L,\infty )$. 
Then $f$ satisfies (M), (S), (C) and (L) with 
 $\pf =p$. Since we are assuming $p>\pt$,  $f$ is also supercritical.  Clearly the functions  $\vp ,\psi\colon [0,\infty )\to [0,\infty )$  satisfy the hypotheses of Theorem~\ref{thm:CS} so,  by \eqref{eq:critsing} and \eqref{eq:ORV1},  there exists  ${\lam}_0 >0$ such that 
\[
\mu_c(x)= \psi\left( \lam_0 |x|^{-\theta} \right)=
\lam_0^{\frac{1}{p-1}} |x|^{-\frac{\theta}{p-1}}\left(\log\left(\lam_0 |x|^{-\theta}\right)\right)^{-\frac{q}{p-1}}\chi_r(x),\qquad x\in\R^N, 
\]
is a DCS, where $r={(\lambda_0/L)^{1/\theta}}$. 
We note that, as in \cite{FHIL}*{Theorem~1.1}, 
\begin{equation}
  \label{eq:1.13}
  \mu_c(x)\asymp
  |x|^{-\frac{\theta}{p-1}}\left|\log |x|\right|^{-\frac{q}{p-1}}\qquad \text{as}\ x\to 0.
\end{equation}
}
\end{Example}

\begin{remark}
\label{rem:nonmonof}
{\rm
Functions of the form $f(u) = u^p(\log (1+u))^q$ with  $p>0$ and $q\in \mathbb{R}$ were also considered in
 \cite{HZ} and \cite{S} as  examples of nonlinearities which are not asymptotically scale-invariant. There   
  the authors were still able to determine  relevant  blow-up rates, universal estimates and Liouville-type theorems. 
 While these nonlinearities satisfy (S), (C) and (L) with $\pf =p$, they do not in general satisfy  (M) when $q<0$.  However, it is straightforward to show 
 that Theorem~\ref{thm:main} remains true for such functions. Indeed, by taking $L>0$ large enough as in Example~\ref{eg:powerlog},  
one may construct functions ${f}_1$ and ${f}_2$ such that ${f}_1(u)\le f(u)\le {f}_2(u)$ for  $u\ge 0$ 
and $f(u)={f}_1(u)={f}_2(u)$ for $u\ge L$, but   
${f}_1$ and ${f}_2$ also satisfy (M) (in addition to (S), (C) and (L)). Then, via  comparison, both parts of Theorem~\ref{thm:main} follow.
}
\end{remark}

\begin{remark}
\label{rem: singorigin}
{\rm 
For purely expositional purposes, we  have chosen to present our results 
with reference to the origin. Of course, by the spatial translational invariance in $\R^N$ of the fractional Laplacian, one can formulate our results with reference to any point $x_0\in\R^N$. For example, defining 
\[
  \CB (x_0) := \bigcap_{r>0} L^{\infty }\left(\R^n\backslash B(x_0,r)\right), 
\]
the corresponding statement regarding existence in Theorem~\ref{thm:main}~(b) (under the same hypotheses on $f$), with  $\mu\in\CB (x_0)\cap L^{1,+}_{{\rm loc}}(\R^n)$, is that there exists an $\varepsilon >0$ such if 
\[
G(\mu(x))\le\varepsilon |x-x_0|^{-\theta} 
\]
for a.a.~$x$ in a neighborhood of $x_0$, then problem~\eqref{eq:SHE} is locally solvable. Likewise,  
under the same hypotheses as Corollary~\ref{cor:CS} one obtains DCS (now defined relative to $\CB (x_0)$) of the form
\[
\mu_c(x)= \psi_f\left(
     \lambda_0|x-x_0|^{-\theta}\right),\qquad x\in\R^N.
     \] 
}
\end{remark}

\vspace{3mm}

The rest of this paper is organized as follows. 
In Section~\ref{sec:solveS} we recall some properties of the semigroup $S(t)$. 
We also introduce the definition of solutions and supersolutions to problem~\eqref{eq:SHE}. 
We obtain necessary conditions and sufficient conditions 
on the  local  solvability of problem~\eqref{eq:SHE} under assumptions~(M), (S) and (C). 
Furthermore, we prepare some lemmas on the inverse function of $F$, 
and we prove Theorem~\ref{thm:main}~(a) using the behavior of $f(F^{-1}(u))$ as $u\to 0$ 
and the convexity of a variant of $F^{-1}$.
Section~\ref{sec:proofs} is devoted to the proof of Theorem~\ref{thm:main} and Theorem~\ref{thm:CS}. Finally, in Section~\ref{section:apps} we consider some applications of our general theory to some specific examples. In particular, we obtain new results for a generalised exponential nonlinearity. 
\section{Local solvability of the  semilinear problem~\eqref{eq:SHE}}
\label{sec:solveS}
In this section we introduce some further notation and terminology,  
defining our  solution concept for  problem~\eqref{eq:SHE} alongside that of a supersolution. 
We recall some basic facts concerning the fractional heat kernel and its associated evolution semigroup.
 Two key theorems  (Theorem~\ref{thm:nec} and Theorem~\ref{thm:suff}) underpin the subsequent proofs of our main results
 regarding necessity and sufficiency for the  local  solvability for problem~\eqref{eq:SHE}. We obtain them by modifying the arguments in~\cites{FI, HI01}, subject to the under conditions (M), (S) and (C). 
\subsection{Preliminaries}
\label{section:2}

For $T>0$, we set $Q_T=\mathbb{R}^N \times (0,T)$.  
We use the abbreviation `a.a.' for `almost all', with respect to Lebesgue measure. 
To avoid unnecessarily cumbersome notation, 
we will use $C$ to denote generic positive constants which are independent of 
 the space variable $x$, and the time variable $t$. In particular, $C$ may take different values within a calculation.
\vspace{5pt}

For any $\theta\in(0,2]$, 
let $\Gamma_\theta=\Gamma_\theta(x,t)$ denote the fundamental solution of the fractional heat equation
\[
\partial_t u+(-\Delta)^{\frac{\theta}{2}}u=0\qquad\mbox{in}\qquad{{\R}}^N\times(0,\infty).
\]
Then  for all $x\in{{\R}}^N$ and $t>0$, we have
\begin{equation}
  \label{eq:kernelbound}
  \begin{split}
    & \Gamma_\theta(x,t)=(4\pi t)^{-\frac{N}{2}}\exp\left(-\frac{|x|^2}{4t}\right)\quad\mbox{if}\quad \theta=2,\\
    & 
    c_\theta t^{-\frac{N}{\theta}}\left(1+t^{-\frac{1}{\theta}}|x|\right)^{-N-\theta}\le
    \Gamma_\theta(x,t)
    \le C_\theta t^{-\frac{N}{\theta}}\left(1+t^{-\frac{1}{\theta}}|x|\right)^{-N-\theta}
    \quad\mbox{if}\quad 0<\theta<2, 
  \end{split}
\end{equation}
where $c_\theta$ and $C_\theta$ are positive constants depending only on $N$ and $\theta$. 
Furthermore, the following properties hold:
\begin{equation}
  \label{eq:2.2}
  \begin{aligned}
    &
    \bullet\quad  \mbox{$\Gamma_\theta$ is positive and smooth in ${{\R}}^N\times(0,\infty)$},
    \\
    &
    \bullet\quad \mbox{$\Gamma_\theta(\cdot,t)$ is radially symmetric}, 
    \\ 
    &
    \bullet\quad   \Gamma_\theta(x,t)=t^{-\frac{N}{\theta}}\Gamma_\theta\left(t^{-\frac{1}{\theta}}x,1\right),
    \qquad \int_{{\mathbb R}^N}\Gamma_\theta(x,t)\,\dee x=1,
    \\
    &
    \bullet\quad  \Gamma_\theta(x,t)=\int_{{{\R}}^N}\Gamma_\theta(x-y,t-s)\Gamma_\theta(y,s)\,\dee y,
  \end{aligned}
\end{equation}
for all $x\in{{\R}}^N$ and $0<s<t$ (see e.g.,  \cites{BJ, BK} for $0<\theta<2$).

\begin{remark}
\rm{
The constant  in the lower bound of \eqref{eq:kernelbound} satisfies $c_\theta\to 0$ as $\theta\to 2$ (see \cite{BG}*{Theorem 2.1}). In particular, the estimate from below in \eqref{eq:kernelbound} is not valid when $\theta=2$ since $\Gamma_2(\cdot,t)$ decays exponentially.  
However,  the estimate from above remains valid for $\theta=2$, by direct comparison with $\Gamma_2$.  Thus, upon relabelling $C_\theta$, we have
\begin{equation}
 0<\Gamma_\theta(x,t)
    \le C t^{-\frac{N}{\theta}}\left(1+t^{-\frac{1}{\theta}}|x|\right)^{-N-\theta}
     \label{eq:2.1}
\end{equation}
for all $x\in{{\R}}^N$, $t>0$ and  $\theta\in (0,2]$. 
This upper bound is sufficient for our purposes. 
We will obtain the relevant lower bounds by different means, via Theorem~\ref{thm:nec}.
} 
\end{remark}

We recall the following lemma on the decay estimate of the semigroup $S(t)$ associated with $\Gamma_\theta$
(see \cite{HI01}*{Lemma~2.1}).
\begin{lemma}
  \label{Lemma:2.1}
  For any  locally integrable function $\mu$ in $\R^N$, 
  set
  \[
  [S(t)\mu](x):= \int_{{{\R}}^N} \Gamma_\theta(x-y,t)\mu(y) \,\dee y, \qquad x\in{{\R}}^N,\,\,\, t>0.
  \]
   If $\mu\in L^1_{{\rm uloc}}\left({\mathbb R}^N\right)$, then there exists $C=C(N,\theta)>0$ such that
  \[
  \|S(t)\mu\|_{L^\infty({{\R}}^N)} \le C t^{-\frac{N}{\theta}} 
  \sup_{x\in{{\R}}^N} \int_{B(x,t^\frac{1}{\theta})} \mu(y)\, \dee y, 
  \qquad t>0.
  \]
\end{lemma}

We now  define our notion of solution for problem~\eqref{eq:SHE}. 
\begin{definition}
  \label{Definition:2.1}
  Let $\mu\in L^{1,+}_{{\rm loc}}(\R^N)$,  $T>0$ and $u$ be a nonnegative, measurable, finite almost everywhere function in $Q_T$.
  Let $f$ be a nonnegative and continuous function in $[0,\infty)$. 
  We say that $u$ is a 
   {\rm solution} of problem~\eqref{eq:SHE} in $Q_T$ if $u$ satisfies 
  \begin{equation}
  \label{eq:2.3}
  u(x,t)=\int_{{{\R}}^N}\Gamma_\theta(x-y,t)\mu(y) \,\dee y+\int_0^t\int_{{{\R}}^N}\Gamma_\theta(x-y,t-s)f(u(y,s))\,\dee y\,\dee s
  \end{equation}
  for a.a.~$(x,t) \in Q_T$. If such a solution exists for some $T>0$, then we say that problem~\eqref{eq:SHE} is {\rm\bf  locally solvable}.

  If $u$ satisfies \eqref{eq:2.3} with `$=$' replaced by `$\ge$',
  then $u$ is said to be a supersolution of problem~\eqref{eq:SHE} in $Q_T$.
\end{definition} 

\vspace{3pt}

The following lemma is crucial in the arguments of this paper.  It asserts that 
the existence of a supersolution of problem~\eqref{eq:SHE} implies the existence of a solution of problem~\eqref{eq:SHE}. 
\begin{lemma}
  \label{Lemma:2.2}
  Let  $\mu\in L^{1,+}_{{\rm loc}}(\R^N)$  and $T>0$. 
  Let $f$ be  a nonnegative, continuous and nondecreasing function in $[0,\infty)$. 
  If there exists a supersolution~$v$ of problem~\eqref{eq:SHE} in $Q_T$,
  then there exists a solution~$u$ of problem~\eqref{eq:SHE} in $Q_T$ such that
  \[
  0\le u(x,t)\le v(x,t),\qquad {\text a.a.\ }(x,t)\in Q_T.
  \]
\end{lemma}
Lemma~\ref{Lemma:2.2} follows from the nonnegativity of  $\Gamma_\theta$ and a standard monotone iteration method.  See for example, \cite{FHIL}*{Lemma~2.2} (also \cite{HI01}*{Lemma~2.2}).

\subsection{A necessary condition for  local solvability}
\label{subsec:solveS.1}

The {\it a priori} esimates we prove below in Lemma~\ref{lem:nec} and Theorem~\ref{thm:nec} are similar in appearance to those previously obtained in the study of semilinear heat equations; see for example  \cite{W3}*{Theorem~1} (or  \cite{QS}*{Lemma~15.6}), \cite{W1}*{Theorem~5}, \cite{Fuj}*{Theorem~1} and  \cite{FI}*{Lemma~4.1}.
Such estimates are essentially based on the monotonicity and convexity of $f$ and positivity of the 
underlying semigroup $S(t)$. However, the existing results require at least one of the following assumptions: 
solutions are classical or time-continuous mild solutions, $f$ is convex in $[0,\infty)$ (facilitating Jensen's inequality), $f$ is Lipschitz continuous (permitting comparison methods from ODE theory) or $f(0)=0$. In contrast,  Theorem~\ref{thm:nec} does not require any of these.   

%
\begin{lemma}
  \label{lem:nec} 
  Let  $\mu\in L^{1,+}_{{\rm loc}}(\R^N)$ and $f\colon [0,\infty)\to [0,\infty)$ be convex and satisfy {\rm (S)}.
    Suppose that  $f$ satisfies:
 \begin{itemize}
  \item[{\rm{(M')}}] $f$ is continuous, nondecreasing and there exists $u_0>0$ such that $f(u)>0$ for all $u\ge  u_0$.  
\end{itemize}
Let $\Sigma_{T}(u_0):=\{(x,t)\in Q_T\, :\, [S(t)\mu](x)\ge u_0  \}$. If problem~\eqref{eq:SHE} is locally solvable 
  in $Q_T$ for some $T>0$,
  then 
  \[
  F\left( [S(t)\mu](x) \right) \ge t
  \]
  for a.a.~$(x,t) \in\Sigma_{T}(u_0)$, where  $F$ is given by \eqref{eq:1.3} for $u \ge u_0$.
 %
\end{lemma}
\textbf{Proof.} 
Suppose  problem~\eqref{eq:SHE} possesses a solution~$u$ in $Q_T$ for some $T>0$.
Let 
\[
E:=\{(x,t)\in Q_T\,:\, \mbox{\eqref{eq:2.3} holds and $u(x,t)<\infty$} \}.
\]
For any $x\in{\mathbb R}^N$, set $E_x:=\{t\in(0,T)\,:\,(x,t)\in E\}$. 
By Definition~\ref{Definition:2.1} and Fubini's theorem 
we see that the set $(0,T)\backslash  E_x$ has Lebesgue measure zero for a.a.~$x\in{\mathbb R}^N$. 

Let $x_*\in{\mathbb R}^N$ be such that the set $(0,T)\backslash  E_{x_*}$ has Lebesgue zero. 
Let $t_*\in E_{x_*}$, and set 
\[
h(t)=[S(t_*)\mu](x_*) + \int_0^t f([S(t_*-s)u(s)](x_*))\, \dee s 
\quad\mbox{for}\quad t\in (0,t_*].
\]
It follows from the  convexity of $f$ and Jensen's inequality  that
\begin{equation*}
\begin{split}
  [S(t_*-t)u(t)](x_*)
  & =[S(t_*)\mu](x_*) + \int_0^t [S(t_*-s)f(u(s))](x_*)\, \dee s\\
  & \ge [S(t_*)\mu](x_*) + \int_0^t f([S(t_*-s)u(s)](x_*))\, \dee s=h(t)
\end{split}
\end{equation*}
for a.a~$t\in E_{x_*}\cap (0,t_*]$. 
Similarly, 
\begin{equation*}
h(t) \le [S(t_*)\mu](x_*) + \int_0^{t_*} [S(t_*-s)f(u(s))](x_*)\, \dee s=u(x_*,t_*)<\infty
\end{equation*}
for $t\in (0,t_*]$. 
Then $h$ is absolutely continuous in $[0,t_*]$ and
\begin{equation}
\label{eq:3.1}
  h'(t) = f([S(t_*-t)u(t)](x_*))\ge f(h(t))
\end{equation}
for a.a~$t\in E_{x_*}\cap[0,t_*]$.

Suppose  
$h(0)=[S(t_*)\mu](x_*)\ge u_0$. It follows from \eqref{eq:3.1} that $h(t)\ge  u_0$, and consequently $f(h(t))>0$  for a.a~$t\in E_{x_*}\cap[0,t_*]$. 
 Hence,
\[
\frac{\dee}{\dee t} F(h(t)) = -\frac{h'(t)}{f(h(t))} \le -1
\]
for a.a.~$t\in E_{x_*}\cap[0,t_*]$.
We deduce that
$0\le F(h(t_*))\le F(h(0)) -t_*$, 
and so 
\[
F([S(t_*)\mu](x_*))=F(h(0))\ge t_*.
\]
Since $t_*\in (0,T)$ and $x_*\in{\mathbb R}^N$ are arbitrary, we obtain 
\begin{equation}
\label{eq:3.2} 
F\left( [S(t)\mu](x) \right) \ge t
\end{equation}
for a.a.~$(x,t) \in \Sigma_{T}(u_0)$, as required. 
Thus Lemma~\ref{lem:nec} follows.
$\hfill\Box$

\medskip 

We now weaken the constraint that $f$ be everywhere convex in Lemma~\ref{lem:nec}, at the expense of strengthening the positivity assumption on $f$. The {\it a priori} estimate is then obtained on a smaller time interval. However, since we are concerned only with local solvability this is not a hinderance.
\begin{theorem}
  \label{thm:nec} 
  Let  $\mu\in L^{1,+}_{{\rm loc}}(\R^N)$ and $T>0$. 
  Suppose $f$ satisfies {\rm (M)}, {\rm (S)} and {\rm (C)}, and let $F$ be as in \eqref{eq:1.3}. 
  Then there exist $C_*=C_*(f)\in (0,1)$ and $T_*=T_*(f)\in (0,T)$ such that 
  if problem~\eqref{eq:SHE} is locally solvable 
  in $Q_T$, then 
  \[
  F\left( [S(t)\mu](x) \right) \ge C_*t
  \]
  for a.a.~$(x,t) \in Q_{T_*}$. 
\end{theorem}
\textbf{Proof.} 
Let $\eta$ be a smooth, nondecreasing  function in $[0,\infty)$ 
such that $0\le\eta\le 1$, $\eta=1$ in $[1,\infty)$
and $\operatorname{supp}\eta \subset (0,\infty)$. 
Set 
\begin{equation*}
\tilde{f}(u)=
\left\{
\begin{array}{ll}
0, &  u\in[0,  \tau_2),
\vspace{7pt}\\
\displaystyle{\int_{  \tau_2}^u \eta(s-  \tau_2)f'(s)\,\dee s}, & u\in[  \tau_2,\infty),
\end{array}
\right.
\end{equation*}
where $\tau_2=\max\{   \tau_0,  \tau_1\}$. By {\rm (S)} and {\rm (C)}, $f'$ is eventually positive and so $\tilde{f}$ satisfies (M') for some  $\tilde{u}_0>0$. Since $0\le\eta\le 1$, we  have
that $\tilde{f}\le f$ in $[0,\infty)$. 
For $u\ge \tau_2+1$, we have  
\begin{equation*}
\tilde{f}(u) = \int_{\tau_2}^{\tau_2+1} \eta(s-\tau_2)f'(s)\, \dee s+ \int_{\tau_2+1}^u f'(s)\, \dee s=  f(u) + \tilde{f}(\tau_2+1) - f(\tau_2+1). 
\end{equation*} 
Clearly, by (M) and (S), $f(u)\to\infty$ as $u\to\infty$ and so there exist $C_*=C_*(f )\in (0,1)$ and $\tau_3\ge\tau_2$ such that  
\begin{equation}
\label{eq: ftildef}
\tilde{f}(u)\ge C_*f(u),\qquad u\ge \tau_3.
\end{equation}
Consequently, $\tilde{f}$ satisfies (S). 
Furthermore, since 
\begin{equation*}
\tilde{f}'(u) = 
\begin{cases}
0, & u\in [0,\tau_2), 
\\ 
\eta(u-\tau_2)f'(u), & u\in [\tau_2,\tau_2+1),
\\ 
f'(u), & u\in [\tau_2+1,\infty), 
\end{cases} 
\end{equation*} 
we see that $\tilde{f}'$ is nondecreasing and so $\tilde{f}$ is convex in $[0,\infty)$. 
Hence  $\tilde{f}$ satisfies the hypotheses of Lemma~\ref{lem:nec} (with $f$ and $u_0$ replaced by   $\tilde{f}$ and  $\tilde{u}_0$, respectively).

Now let $\tilde{F}(u)$ be defined by \eqref{eq:1.3} for $u>\tilde{u}_0$, with $f$ replaced by $\tilde{f}$ .  By \eqref{eq: ftildef} we  have
\begin{equation}
  \label{eq:FtildeF}
F(u)\ge C_*\tilde{F}(u),\qquad u\ge\tau_4:=\max\{\tau_3, \tilde{u}_0 \}.
\end{equation}
Set $T_*=T_*(f):=\min\{ T, F(\tau_4)/C_*\}>0$ and 
\[
{\Sigma}_{T_*}(\tau_4):=\{(x,t)\in Q_{T_*}\, :\, [S(t)\mu](x)\ge \tau_4  \}.
\]
Then we have
\begin{equation}
 \label{eq:musmall}
F\left( [S(t)\mu](x) \right)\ge F(\tau_4)\ge C_*T_*
\ge C_*t,\qquad (x,t) \in Q_{T_*}\backslash {\Sigma}_{T_*}(\tau_4).
\end{equation}

Since $\tilde{f}\le f$ and problem~\eqref{eq:SHE} possesses a solution~$u$ in $Q_T$ for some $T>0$,  
Lemma~\ref{Lemma:2.2} implies the existence of a solution $\tilde{u}$ of problem~\eqref{eq:SHE} in $Q_{T_*}$, with $f$ replaced by $\tilde{f}$. Hence,
 by \eqref{eq:FtildeF} and Lemma~\ref{lem:nec}, 
 \begin{equation}
  \label{eq:mubig}
F\left( [S(t)\mu](x) \right)\ge C_*\tilde{F}\left( [S(t)\mu](x) \right)\ge C_*t,\qquad a.a.\ (x,t) \in {\Sigma}_{T_*}(\tau_4).
\end{equation} 
Theorem~\ref{thm:nec} now follows from \eqref{eq:musmall} and \eqref{eq:mubig}.
$\hfill\Box$



\subsection{A sufficient condition for  local solvability}
\label{subsec:solveS.2}

We now obtain sufficient conditions for the  local solvability of problem~\eqref{eq:SHE}.
Thanks to Lemma~\ref{Lemma:2.2}, 
it suffices to find a sufficient condition for the existence of a supersolution of problem~\eqref{eq:SHE}. 
Henceforth we denote by ${\dashint_{B}g}$  the average value of any real-valued function $g$ over a ball $B\subset\R^N$. 
\begin{Theorem}
  \label{thm:suff} 
  Let $f$ satisfy {\rm (M)} and {\rm (S)} and  let  $\mu\in L^{1,+}_{{\rm loc}}(\R^N)$.
  Suppose there exist $\beta>0$, $\delta>0$ and $ \tau^*\in(  \tau_0,\infty)$ such that 
  \begin{equation}
    \label{eq:3.4}
    \beta \le f'(u)F(u) \le 1+\beta-\delta
  \end{equation}
  for all $u\ge  \tau^*$.
  Then, for any $\kappa >1$, there exist $\varepsilon>0$ and $T_*>0$ such that if $\mu$  satisfies
  \begin{equation}  
    \label{eq:3.5} \sup_{x\in\R^N} \dashint_{B(x,\sigma)} G(\mu(y))^\beta \, \dee y 
    \le  
    \varepsilon\sigma^{-\beta\theta},\qquad 0<\sigma<T^{\frac{1}{\theta}}, 
  \end{equation}
  for some $T\in(0,T_*)$, then problem~\eqref{eq:SHE} possesses a solution $u$ in $Q_T$ such that
  \[
  0\le G(u(x,t))^\beta\le\kappa [S(t)\left(G(\mu)^\beta\right)](x),\qquad  {\text a.a.\ }(x,t)\in Q_T. 
  \]
Furthermore, $u\in L_{\rm{loc}}^\infty\left((0,T),L^\infty (\R^N)\right)$.
\end{Theorem}
\textbf{Proof.} 
We can assume, without loss of generality, that $\mu(x)\ge  \tau^*$ for $x\in \R^N$.
Indeed, by Lemma~\ref{Lemma:2.2} it suffices to prove the existence of a solution of  problem~\eqref{eq:SHE}
with $\mu$ replaced by 
$\overline{\mu}:=\max\{\mu, \tau^*\}$. 
Assume that \eqref{eq:3.5} holds for some $\varepsilon>0$ and $T>0$.  
Then, 
since $G(\overline{\mu})^\beta \le G(\tau^*)^\beta+G(\mu)^\beta$ in $\mathbb{R}^N$, 
taking a sufficiently small $T>0$ if necessary, by \eqref{eq:3.5} we have 
\begin{align*}
  \sup_{x\in\R^N} \dashint_{B(x,\sigma)} G(\overline{\mu}(y))^\beta \, \dee y 
  & 
  \le 
  \sup_{x\in\R^N} \dashint_{B(x,\sigma)} G(\tau^*)^\beta \, \dee y 
  + 
  \sup_{x\in\R^N} \dashint_{B(x,\sigma)} G(\mu(y))^\beta \, \dee y 
  \\ 
  & 
  \le G(\tau^*)^\beta  
  + 
  \varepsilon \sigma^{-\beta\theta}
  \le 
  2\varepsilon \sigma^{-\beta\theta}
\end{align*}
for all $\sigma\in (0,T^\frac{1}{\theta})$.  
Thus, the estimate \eqref{eq:3.5} is still valid for $\overline{\mu}$ with $\varepsilon$ replaced by $2\varepsilon$. 

Set 
  \begin{equation}
    \label{eq:3.6}
    \Phi(u)=G(u)^\beta\qquad\mbox{for $u\in[\tau^*,\infty)$}. 
  \end{equation}
It easily follows from \eqref{eq:1.3} that 
$\Phi$ is strictly increasing in $(\tau^*,\infty)$ and $\Phi(u)\to\infty$ as  $u\to\infty$. 
Furthermore, by \eqref{eq:3.4} we see that
\begin{equation*}
  \begin{split}
    \Phi''(u) 
    & =(F(u)^{-\beta})''= \left(\beta F(u)^{-\beta-1}f(u)^{-1}\right)'\\
    & = \beta(\beta+1)F(u)^{-\beta-2}f(u)^{-2}-\beta F(u)^{-\beta-1}f(u)^{-2}f'(u)\\
    & = \frac{\beta \left( (\beta + 1) - f'(u)F(u) \right)}{f(u)^2 F(u)^{\beta+2}}
    >0
  \end{split}
\end{equation*}
for $u\ge  \tau^*$, so that 
\begin{equation}
 \label{eq:3.7}
 \mbox{$\Phi$ is convex in $( \tau^*,\infty)$.}
\end{equation}

Let $\varepsilon>0$ and $T_*>0$ be small, to be chosen later. 
By \eqref{eq:3.5} we apply Lemma~\ref{Lemma:2.1} to obtain 
\begin{equation}
  \label{eq:3.8}
  \| S(t)\Phi(\mu) \|_\infty
  \le
  C t^{-\frac{N}{\theta}}
  \sup_{x\in{{\R}}^N}\int_{B(x,t^{1/\theta})}
  \Phi(\mu(y)) \, \dee y
  \le C\varepsilon t^{-\beta}
\end{equation}
for $t\in(0,T)$, where $T$ is as in \eqref{eq:3.5}.  
Since
\[
[S(t)\Phi(\mu)](x)
\ge [S(t)\Phi( \tau^*)](x) 
= \Phi( \tau^*)>\Phi(0)\ge 0,\qquad(x,t)\in {Q_T},
\]
for any $\kappa >1$, we can define the nonnegative function $v$ in {$Q_T$ }
by
\begin{equation}
\label{eq:3.9}
v(x,t)=\Phi^{-1}\left(\kappa [S(t)\Phi(\mu)](x) \right),
\qquad
(x,t)\in Q_T.
\end{equation}
Notice that the domain of $\Phi^{-1}$ is $(\Phi(0),\infty)$. 
Since $\Phi^{-1}$ is strictly increasing,
by Jensen's inequality, \eqref{eq:3.7} and \eqref{eq:3.9} we have
\begin{equation}
  \label{eq:3.10}
  \begin{split}
    \Phi([S(t)\mu](x))
     & \le
    [S(t)\Phi(\mu)](x),\quad\text{and}\\
    {[S(t)\mu](x)}
     & \le \Phi^{-1}\left([S(t)\Phi(\mu)](x)\right)
    \le
    \Phi^{-1}\left(\kappa [S(t)\Phi(\mu)](x)\right)
    =
    v(x,t),
  \end{split}
\end{equation}
for $(x,t)\in {Q_T }$. 
Since 
$\Phi^{-1}(u)=F^{-1}( u^{-1/\beta})$ for $u>0$ (see \eqref{eq:1.3} and \eqref{eq:3.6}),
by \eqref{eq:3.8} and \eqref{eq:3.9} we obtain
\begin{equation}
  \label{eq:3.11}
  \|v(t)\|_\infty\le 
  \Phi^{-1}\left(\kappa C\varepsilon t^{-\beta}\right)
  =
  F^{-1}
  \left(\kappa^{-\frac{1}{\beta}}C^{-\frac{1}{\beta}}\varepsilon^{-\frac{1}{\beta}} t\right)
  =:\xi(t)
\end{equation}
for $t\in(0,T)$. 
On the other hand,
it follows from \eqref{eq:3.7} that
\begin{equation}
  \label{eq:3.12}
  \mbox{$\Phi'(u)=\beta f(u)^{-1}F(u)^{-(\beta+1)}$\quad
  is strictly increasing in $( \tau^*,\infty)$}.
\end{equation}
Then, by \eqref{eq:3.10}, \eqref{eq:3.11} and \eqref{eq:3.12}
we apply the mean value theorem to obtain
\begin{equation}
  \label{eq:3.13}
  \begin{aligned}
    & (\kappa - 1) [S(t)\Phi(\mu)](x)
    \le 
    \Phi(v(x,t)) - \Phi([S(t)\mu](x))\\
    & \qquad\qquad
    \le
    \Phi'(\xi(t))(v(x,t)-[S(t)\mu](x))
    =
    \frac{\beta\left( v(x,t)-[S(t)\mu](x)\right)}{f(\xi(t))F(\xi(t))^{\beta+1}}
    \\
    & \qquad\qquad
    =C \beta \kappa\varepsilon
    \frac{\left( v(x,t)-[S(t)\mu](x)\right)}{t^{\beta}f(\xi(t))F(\xi(t))}
  \end{aligned}
\end{equation}
for $(x,t)\in {Q_T }$. Now observe from \eqref{eq:3.4} that
\begin{equation*}
  \frac{\dee}{\dee u} \left[f(u)F(u)^\beta\right]
  = F(u)^{\beta-1}\left[ f'(u)F(u) - \beta\right]
  \ge 0
\end{equation*}
for $u\ge  \tau^*$.
Then, by \eqref{eq:3.9} and \eqref{eq:3.11} we see that
\begin{equation*}
\begin{split}
  f(v(x,s))
  & =
  \kappa f(v(x,s)) 
   \left\{ \kappa [S(s)\Phi(\mu)](x) \right\}^{-1}[S(s)\Phi(\mu)](x)
  \\
  &
  =\kappa  
  f(v(x,s)) F(v(x,s))^\beta 
   [S(s)\Phi(\mu)](x)\\
  & 
  \le\kappa 
  f(\xi(s))F(\xi(s))^{\beta-1}F(\xi(s))[S(s)\Phi(\mu)](x)\\
  & 
  =
  C\kappa^{\frac{1}{\beta}}\varepsilon^{\frac{1-\beta}{\beta}}
  s^{\beta-1}f(\xi(s))F(\xi(s))[S(s)\Phi(\mu)](x)
\end{split}
\end{equation*}
for $(x,s)\in Q_T$.
Hence,
\begin{equation}
  \label{eq:3.14}
  \int_0^t [S(t-s)f(v(s))](x) \, \dee s
  \le C\kappa^{\frac{1}{\beta}}\varepsilon^{\frac{1-\beta}{\beta}} [S(t)\Phi(\mu)](x)
  \int_0^t s^{\beta-1} f(\xi(s))F(\xi(s)) \, \dee s
\end{equation}
for $t\in(0,T)$.
Since
\begin{equation*}
  {\left(F^{-1}(u)\right)'}
  = - f(F^{-1}(u)),
\end{equation*}
taking small enough $T_*>0$ if necessary, 
by \eqref{eq:3.4} and \eqref{eq:3.11} we have
\begin{equation}
  \notag
  \begin{aligned}
    &
    \left\{ s^{\beta-\delta} f(\xi(s))F(\xi(s)) \right\}'\\
    & =
    (\beta-\delta)s^{\beta-\delta-1}f(\xi(s))F(\xi(s))
    + s^{\beta-\delta} \left[ f'(\xi(s))F(\xi(s)) - 1 \right] \xi'(s)
    \\
    & =
    (\beta-\delta)s^{\beta-\delta-1}f(\xi(s))
    C^{-\frac{1}{\beta}}\varepsilon^{-\frac{1}{\beta}} \kappa^{-\frac{1}{\beta}} s
    - s^{\beta-\delta} \left[ f'(\xi(s))F(\xi(s)) - 1 \right]
    C^{-\frac{1}{\beta}}\varepsilon^{-\frac{1}{\beta}} \kappa^{-\frac{1}{\beta}} f(\xi(s))
    \\
    & = C^{-\frac{1}{\beta}}\varepsilon^{-\frac{1}{\beta}} \kappa^{-\frac{1}{\beta}}
    s^{\beta-\delta} f(\xi(s)) \left[ (\beta-\delta)
    - \left\{ f'(\xi(s))F(\xi(s)) - 1 \right\} \right]
    \ge 0
  \end{aligned}
\end{equation}
for $s\in(0,T)\subset(0,T_*)$. 
This, together with \eqref{eq:3.14}, implies that
\begin{equation}
  \notag 
  \begin{aligned}
    \int_0^t [S(t-s)f(v(s))](x) \, \dee s
    & \le
    C\kappa^{\frac{1}{\beta}}\varepsilon^{\frac{1-\beta}{\beta}} [S(t)\Phi(\mu)](x)
     t^{\beta-\delta} f(\xi(t)) F(\xi(t))
    \int_0^t s^{\delta-1} \, \dee s
    \\
    & \le 
    C\kappa^{\frac{1}{\beta}}\varepsilon^\frac{1}{\beta}
     \varepsilon^{-1} t^\beta f(\xi(t))F(\xi(t))
    [S(t)\Phi(\mu)](x)
  \end{aligned}
\end{equation}
for $t\in(0,T)$.
Combining this with \eqref{eq:3.13} and taking small enough $\varepsilon\in(0,1)$ if necessary, 
we obtain
\begin{equation}
  \notag
  \begin{aligned}
    v(x,t)-[S(t)\mu](x)
    & \ge
    C^{-1} \beta^{-1} (1-\kappa^{-1}) 
     \varepsilon^{-1}
    t^\beta f(\xi(t))F(\xi(t))
    [S(t)\Phi(\mu)](x)
    \\
    & \ge
    C^{-1} \beta^{-1} (1-\kappa^{-1}) 
    C^{-1} \varepsilon^{-\frac{1}{\beta}} \kappa^{-\frac{1}{\beta}}
    \int_0^t [S(t-s)f(v(s))](x) \, \dee s\\
    & \ge
    \int_0^t [S(t-s)f(v(s))](x) \, \dee s
  \end{aligned}
\end{equation}
for $(x,t)\in Q_T$. 
Consequently $v$ is a supersolution of problem~\eqref{eq:SHE} in $Q_T$. 
Hence, by Lemma~\ref{Lemma:2.2}, we obtain a solution $u$ of problem~\eqref{eq:SHE} in $Q_T$ such that
\[
  0\le u(x,t)\le v(x,t)=\Phi^{-1}\left(\kappa [S(t)\Phi(\mu)](x) \right),
  \qquad (x,t)\in Q_T.
\]
By Lemma~\ref{Lemma:2.1},   $u\in L_{\rm{loc}}^\infty\left((0,T),L^\infty (\R^N)\right)$   
and  Theorem~\ref{thm:suff} follows.
$\hfill\Box$
\section{Proofs of the main results}
\label{sec:proofs}


Assume (M) and (S) and let $F$ be as in \eqref{eq:1.3}.   We first record some basic properties of $F$.   
Since $F$ is strictly decreasing in $(0,\infty)$ and $F(u)\to 0$ as $u\to\infty$, 
the inverse function $F^{-1}$ of $F$ satisfies the following properties: 
\begin{align}
\label{eq:4.1}
 & \text{$F^{-1}$ is strictly decreasing in $(0,F_0)$};\\
\label{eq:4.2}
& \text{$F^{-1}(\sigma)\to\infty$ as $\sigma\to +0$};\\
\label{eq:4.3}
&  \left\{F^{-1}(\sigma)\right\}' = -f\left(F^{-1}(\sigma)\right)\mbox{ for $\sigma \in (0,F_0)$}.
\end{align}
%

\subsection{Preparatory lemmas}
\label{subsection:lemmas}


Here we determine the behavior of $f(F^{-1}(\sigma))$ and $F^{-1}(\sigma)$ as $\sigma\to +0$ 
and obtain the convexity of the function $F^{-1}(\sigma^k)$, with $k>0$. 
\begin{lemma}
\label{Lemma:4.1}
Suppose  $f$ satisfies {\rm (M)}, {\rm (S)} and {\rm (L)}. 
Then, for any $\varepsilon\in(0,q_f)$, there exists $C\ge 1$ such that 
\begin{equation}
\label{eq:4.4}
  C^{-1}\sigma^{\varepsilon-{\qf}}\le f\left(F^{-1}\left(\sigma\right)\right)\le C\sigma^{-\varepsilon -{\qf}}
  \end{equation}
  and
  \begin{equation}
\label{eq:4.5}
F^{-1}(\sigma)\le C\sigma^{1-{\qf}-\varepsilon},
\end{equation}
for small enough $\sigma>0$. 
\end{lemma}
{\bf Proof.}
Set $h(\sigma):=f\left(F^{-1}(\sigma)\right)$ for $\sigma\in(0,F_0)$. 
It follows from \eqref{eq:4.3} that 
\begin{equation}
\label{eq:4.6}
\left\{F^{-1}(\sigma)\right\}'=-h(\sigma)\mbox{ for $\sigma \in (0,F_0)$}.
\end{equation}
Hence,
\begin{equation}
\label{eq:4.7}
\begin{split}
h'(\sigma) & =f'\left(F^{-1}(\sigma)\right)\left\{F^{-1}(\sigma)\right\}'=-f'\left(F^{-1}(\sigma)\right)h(\sigma)
=-\frac{f'(u)F(u)}{F(u)}h(\sigma)\\
 & =-f'(u)F(u)\sigma^{-1}h(\sigma)
\end{split}
\end{equation}
for $\sigma\in(0,F_0)$ with $u=F^{-1}(\sigma)$. 

Let $\varepsilon>0$. 
It follows from ${\qf}\ge 1$ (see Remark~\ref{rem:ORV}) that
\[
{\qf}-\varepsilon\le f'(u)F(u)\le {\qf}+\varepsilon
\] 
for large enough $u>0$. 
Then, by \eqref{eq:4.2} and \eqref{eq:4.7} we find $\sigma_*>0$ such that 
\[
-({\qf}+\varepsilon)\sigma^{-1}h(\sigma)\le h'(\sigma)\le -({\qf}-\varepsilon)\sigma^{-1}h(\sigma)
\]
for $\sigma\in(0,\sigma_*)$. 
Then we observe that
\begin{align*}
\left(\sigma^{{\qf}+\varepsilon}h(\sigma)\right)' & =({\qf}+\varepsilon)\sigma^{{\qf}+\varepsilon-1}h(\sigma)+\sigma^{{\qf}+\varepsilon}h'(\sigma)\ge 0,\\
\left(\sigma^{{\qf}-\varepsilon}h(\sigma)\right)' & =({\qf}-\varepsilon)\sigma^{{\qf}-\varepsilon-1}h(\sigma)+\sigma^{{\qf}-\varepsilon}h'(\sigma)\le 0,
\end{align*}
for $\sigma\in(0,\sigma_*)$. 
These imply that
\begin{equation}
\label{eq:4.8}
\sigma^{{\qf}+\varepsilon}h(\sigma)\le \sigma_*^{{\qf}+\varepsilon}h(\sigma_*)
\qquad\text{and}\qquad 
\sigma^{{\qf}-\varepsilon}h(\sigma)\ge \sigma_*^{{\qf}-\varepsilon}h(\sigma_*),
\end{equation}
for $\sigma\in(0,\sigma_*)$. Thus \eqref{eq:4.4} follows. 

Finally, since ${\qf}\ge 1$, 
taking small enough $\sigma_*>0$ if necessary, 
by \eqref{eq:4.6} and \eqref{eq:4.8} we obtain
\begin{align*}
F^{-1}(\sigma) & =F^{-1}(\sigma_*)-\int_{\sigma_*}^\sigma h(s)\, \dee s=F^{-1}(\sigma_*)+\int^{\sigma_*}_\sigma h(s)\,\dee s\\
 & \le C+C\int^{\sigma_*}_\sigma s^{-{\qf}-\varepsilon}\,\dee s
\le C\sigma^{-{\qf}+1-\varepsilon}
\end{align*}
for $\sigma\in(0,\sigma_*)$. Thus \eqref{eq:4.5} follows. 
$\hfill\Box$
\begin{lemma}
  \label{Lemma:4.2}
  Suppose the hypotheses of Lemma~{\rm\ref{Lemma:4.1}} hold. 
  For $k>0$, set 
  \begin{equation}
    \label{eq:4.9}
    g(\sigma) = F^{-1}\left(\sigma^k\right),\qquad \sigma\in\left(0,F_0^{\frac{1}{k}}\right).
  \end{equation}
  Then there exists $\sigma_*>0$ such that 
  $g$ is convex in $(0,\sigma_*)$. 
\end{lemma}

\noindent 
\textbf{Proof.}
It follows from Remark~\ref{rem:ORV} that ${\qf} \ge 1$. 
Then, for any $k>0$, we have 
\begin{equation}
  \label{eq:4.10}
  f'(u)F(u) \ge 1-\frac{1}{2k} 
\end{equation}
for large enough $u>0$. 
By \eqref{eq:4.3} we see that 
\begin{align*}
  g'(\sigma)
  &=
  -k f\left(F^{-1}\left(\sigma^k\right)\right)
  \sigma^{k-1}, 
  \\ 
  g''(\sigma)
  &= 
  k^2 
  f'\left(F^{-1}\left(\sigma^k\right)\right)
  f\left(F^{-1}\left(\sigma^k\right)\right)
  (\sigma^{k-1})^2
  - k(k-1)
  f\left(F^{-1}\left(\sigma^k\right)\right)\sigma^{k-2}
  \\ 
  &= 
  k^2 f\left(F^{-1}\left(\sigma^k\right)\right)
  \sigma^{k-2} \left[ f'\left(F^{-1}\left(\sigma^k\right)\right)
  F\left(F^{-1}\left(\sigma^k\right)\right) - \left(1-\frac{1}{k}\right) \right],
\end{align*}
for $\sigma\in(0,F_0^{\frac{1}{k}})$. 
Then, by \eqref{eq:4.2} and \eqref{eq:4.10} we find $\sigma_*>0$ such that $g''(\sigma)>0$ for $\sigma \in (0,\sigma_*)$. 
$\hfill\Box$\vspace{5pt}

\subsection{Proof of Theorem~\ref{thm:main}}
\label{subsection:pf1}


We first prove part (a) of Theorem~\ref{thm:main} (nonexistence) under the assumptions {\rm (M)}, {\rm (S)}, {\rm (C)} and {\rm (L)}.
Let   $\gamma\ge 1$ (to be chosen later) and suppose $\mu$ satisfies
\[
G(\mu(x))\ge\gamma |x|^{-\theta},\qquad x\in B_R,
\]
for some  $R>0$. 
Suppose also that problem~\eqref{eq:SHE} is locally solvable.
Let $r\in (0,R)$ be such that $r^{-\theta}> G_0$ (so that ${\gamma} |x|^{-\theta}>G_0$ for $|x|<r$) and set 
\begin{equation*}
  \hat{\mu}(x) = 
  \Ginvext \left( {\gamma} |x|^{-\theta} \right) \chi_r(x), 
  \qquad x\in\mathbb{R}^N, 
\end{equation*}
with $\Ginvext$ as in Remark~\ref{rem:ORV}~(iv). Clearly $\hat{\mu}\le\mu$  for $|x|\ge r$ and  
\begin{equation*}
  \hat{\mu}(x) = G^{-1} \left( {\gamma} |x|^{-\theta} \right),\qquad x\in B_r, 
\end{equation*}
so that $\hat{\mu}\le\mu$ in ${\mathbb R}^N$. 
Hence, by Theorem~\ref{thm:nec},  there exist $C_*=C_*(f)>0$ and $T_*=T_*(f)>0$ such that
\begin{equation}
  \label{eq:4.11}
  F(\| S(t)\hat{\mu}\|_\infty)\ge  F (\| S(t)\mu\|_\infty)\ge C_*t 
\end{equation}
for a.a.~$t\in (0,T_*)$.

Now fix $\alpha \in (0,\theta)$ and let $g$ be as in \eqref{eq:4.9} with $k=\theta/\alpha >1$; that is, 
\[
g(\sigma)=F^{-1}\left(\sigma^{\frac{\theta}{\alpha}}\right),\qquad \sigma\in (0,\sigma_*),
\]
where $\sigma_*>0$ is chosen such that (by Lemma~\ref{Lemma:4.2})  $g$ is convex in $(0,\sigma_*)$. 
In addition we choose $r, \sigma_*>0$ small enough such $\sigma_*^{-1}>G_0$ and $r<{\sigma_*}^{1/\alpha}$. 
Note in particular that for $|x|<r$, we have  $\gamma^{-\alpha /\theta}|x|^{\alpha} <\sigma_*$.  
Since 
\[
G^{-1}\left(\sigma^{-1}\right) = F^{-1}(\sigma) = g\left(\sigma^{\frac{\alpha}{\theta}}\right),\qquad 
 \sigma\in\left(0,\min\{G_0^{-1}, \sigma_*^{\theta/\alpha }\}\right),  
\]
we have 
\[
\hat{\mu}(x) =
  G^{-1}\left(\left(\gamma^{-1}|x|^\theta\right)^{-1}\right)=
  g \left(\gamma^{-\frac{\alpha}{\theta}}|x|^\alpha\right),
   \qquad x\in B_r. 
\]
Consequently,
\[
\hat{\mu}(x) =
\left\{
\begin{array}{ll}
 g \left(\gamma^{-\frac{\alpha}{\theta}}|x|^\alpha\right), & |x|<r ,\\
0, & |x|\ge r.
\end{array}
\right.  
\]
This together with \eqref{eq:2.2} implies that, for all $t>0$,
\begin{equation}
\label{eq:4.12}
  \begin{split}
    \| S(t)\hat{\mu} \|_\infty
    \ge 
    [S(t)\hat{\mu}](0) 
    & = 
    \int_{B_r} \Gamma_\theta(y,t)  g \left(\gamma^{-\frac{\alpha}{\theta}} |y|^\alpha\right) \, \dee y 
    \\ 
    & = 
    \int_{B(0,rt^{-\frac{1}{\theta}})} \Gamma_\theta(z,1)  g \left(\gamma^{-\frac{\alpha}{\theta}} t^\frac{\alpha}{\theta} |z|^\alpha\right) \, \dee z.
  \end{split}
\end{equation}
Now set 
\begin{equation}
  \notag 
  m(t) = \int_{B(0,rt^{-\frac{1}{\theta}})} \Gamma_\theta (z,1)\, \dee z,\qquad t>0. 
\end{equation}
Since  $r<\sigma_{*}^{1/\alpha}$ we have that  
$
  \gamma^{-\frac{\alpha}{\theta}} t^\frac{\alpha}{\theta} |z|^\alpha
  <\sigma_*$
  {for all} $ z\in B(0,rt^{-\frac{1}{\theta}})$. 
Hence, by the convexity of $ g $ in $(0,\sigma_*)$ and \eqref{eq:4.12}, we may apply Jensen's inequality to obtain 
\begin{equation}
\label{eq:4.13}
\| S(t)\hat{\mu} \|_\infty
\ge  m(t) 
 g \left(\frac{1}{m(t)} \int_{B(0,rt^{-\frac{1}{\theta}})} \Gamma_\theta(z,1) 
\gamma^{-\frac{\alpha}{\theta}} t^\frac{\alpha}{\theta} |z|^\alpha \, \dee z \right).
\end{equation}
On the other hand, by  \eqref{eq:2.1} and \eqref{eq:2.2} 
we have, for any $t>0$, 
\begin{eqnarray}
\label{eq:4.14}
0\le 1-m(t)=\int_{{\mathbb R}^N\backslash  B(0,rt^{-\frac{1}{\theta}})} \Gamma_\theta (z,1)\, \dee z
&\le& C\int_{{\mathbb R}^N\backslash  B(0,rt^{-\frac{1}{\theta}})}(1+|z|)^{-N-\theta}\, \dee z\nonumber\\
&\le& C_1 r^{-\theta}t
\end{eqnarray}
Choosing   small enough  $t$ (i.e.,  such that $C_1
r^{-\theta}t\le 1/2$), we can ensure that $m(t)\ge 1/2$.   Noting that 
 $\alpha\in (0,\theta)$, 
by \eqref{eq:2.1}  we find $C_2>0$ such that
\begin{equation}
\label{eq:4.15}
\frac{1}{m(t)} \int_{B(0,rt^{-\frac{1}{\theta}})} 
\Gamma_\theta(z,1) \gamma^{-\frac{\alpha}{\theta}} t^\frac{\alpha}{\theta} |z|^\alpha \, \dee z
\le 2 \gamma^{-\frac{\alpha}{\theta}} t^\frac{\alpha}{\theta}
\int_{{\mathbb R}^N}\Gamma_\theta(z,1)|z|^\alpha \, \dee z
\le C_2\gamma^{-\frac{\alpha}{\theta}} t^\frac{\alpha}{\theta}.
\end{equation}
Then, since $ g $ is  non-increasing, we observe from \eqref{eq:4.13}, \eqref{eq:4.14} and \eqref{eq:4.15} that 
\begin{equation}
\label{eq:4.16}
\| S(t)\hat{\mu} \|_\infty
\ge (1-C_1't) g \left(C_2\gamma^{-\frac{\alpha}{\theta}} t^\frac{\alpha}{\theta}\right)
    \ge 
    {F}^{-1}\left(C_3\gamma^{-1}t\right) 
    -C_1' t{F}^{-1}\left(C_3\gamma^{-1}t\right) 
\end{equation}
for  small enough $t>0$, where $C_1'=C_1 r^{-\theta}$.
On the other hand, by \eqref{eq:4.1} and \eqref{eq:4.11} we have 
\begin{equation}
  \notag 
  \| S(t)\hat{\mu} \|_\infty 
  = 
  F^{-1}\left(
    F(\| S(t)\hat{\mu} \|_\infty)
  \right)
  \le 
  F^{-1}(C_*t)
\end{equation}
for  small enough $t>0$. 
Then we deduce from \eqref{eq:4.16} that
\begin{equation}
  \label{eq:4.17}
  F^{-1}\left(C_3\gamma^{-1}t\right) - F^{-1}(C_*t) 
  \le 
  C_1' t F^{-1}\left(C_3\gamma^{-1}t\right)
\end{equation}
for  small enough $t>0$. 
In particular,  we note that the constants $C_*$, $C_1'$ and $C_3$ in \eqref{eq:4.17} are independent of $\gamma$, while $C_*$ and $C_3$  are independent of $r$.

Now let $\varepsilon'\in (0,1/2)$ and choose  $\gamma\ge 2C_3/C_*$ (independent of $r$).
By~\eqref{eq:4.3} we apply the mean value theorem to obtain 
\[
  F^{-1}\left(C_3\gamma^{-1}t\right) - F^{-1}(C_*t) 
  =(C_*-C_3\gamma^{-1})t f\left(F^{-1}(t_*)\right) 
  \ge\frac{C_*}{2} tf\left(F^{-1}(t)\right), 
\]
where $t_*\in(C_3\gamma^{-1}t,C_*t)$. 
This together with \eqref{eq:4.4} implies that 
\begin{equation}
  \label{eq:4.18}
  F^{-1}\left(C_3\gamma^{-1}t\right) - F^{-1}(C_*t) 
  \ge C t^{1-{\qf}+\varepsilon'}
\end{equation}
for  small enough $t>0$. 
On the other hand, 
by \eqref{eq:4.5} we have 
\begin{equation}
  \label{eq:4.19}
  C_1't F^{-1}\left(C_3\gamma^{-1}t\right)
  \le 
  C \gamma^{{\qf}-1+\varepsilon'}t^{2-{\qf}-\varepsilon'}
\end{equation}
for  small enough $t>0$. 
Combining \eqref{eq:4.17}, \eqref{eq:4.18} and \eqref{eq:4.19}, 
we obtain 
\begin{equation}
  \notag 
  0<C\le\gamma^{{\qf}-1+\varepsilon'}
  t^{1-2\varepsilon'}
\end{equation}
for  small enough $t>0$. 
Since $\varepsilon'\in (0,1/2)$, we obtain a contradiction for  small enough $t>0$. 
Thus problem~\eqref{eq:SHE} is not locally solvable.
\vspace{5pt}

Now we apply Theorem~\ref{thm:suff} and Lemma~\ref{Lemma:2.2} 
to prove part (b) of Theorem~\ref{thm:main} (existence).
%
%
Suppose  \eqref{eq:1.5} holds in some neighbourhood of the origin (and therefore in $B_R$, for some $R>0$)  for some $\varepsilon >0$  to be chosen.  
Since $\mu\in\CB$, there exists $K=K(R)\ge G_0$ such that $\mu(x)\le G^{-1}(K)$ for all $x\in\R^N\backslash  B_R$. 

Set
\begin{equation}
  \label{eq:4.20}
  \tilde{\mu}(x) = G^{-1}\left(\varepsilon |x|^{-\theta} + K\right), 
  \qquad x\in\mathbb{R}^N. 
\end{equation}
It follows from \eqref{eq:1.5} and \eqref{eq:4.20} that $\mu\le \tilde{\mu}$ in ${\mathbb R}^N$. 
By Lemma~\ref{Lemma:2.2} it suffices to prove the existence of
a local-in-time solution of problem~\eqref{eq:SHE} with the initial 
function $\tilde{\mu}$. 

Since $f$ is supercritical, $q_f<1+{N}/{\theta}$ and so we may choose $\beta >0$
such that 
\[
{\qf}-1<\beta<\min \{{\qf}, N/\theta  \}.
\]
We may then choose small enough $\delta>0$ so that 
\[
{\qf}-1<\beta<{\qf}<1+\beta-\delta .
\]
Then \eqref{eq:3.4} holds for  large enough $u>0$. 
Furthermore, by \eqref{eq:3.6} and \eqref{eq:4.20} 
we have
\begin{equation*}
  \Phi\left(\tilde{\mu}(x)\right) 
  = 
  G(\tilde{\mu}(x))^{\beta}=\left(\varepsilon |x|^{-\theta} + K\right)^\beta
  \le 
  C\varepsilon^\beta |x|^{-\beta\theta} + CK^\beta
\end{equation*}
for a.a.~$x\in{\mathbb R}^N$ (where $C=C(\beta)$). 
Then, since $\beta\theta<N$, 
we see that
\begin{equation*}
  \sup_{x\in\mathbb{R}^N} \dashint_{B(x,\sigma)} \Phi(\tilde{\mu}(y)) \, \dee y 
  \le 
  C\varepsilon^\beta \sigma^{-\beta\theta} + CK^\beta
  \le C\varepsilon^\beta \sigma^{-\beta\theta}
\end{equation*}
for  small enough $\sigma>0$. 
Taking  small enough $\varepsilon>0$ 
and applying Theorem~\ref{thm:suff}, we see that 
problem~\eqref{eq:SHE}  possesses a solution $\tilde{u}\in L_{\rm{loc}}^\infty\left((0,T),L^\infty (\R^N)\right)$ 
for the initial function $\tilde{\mu}$; 
consequently problem~\eqref{eq:SHE} possesses a solution $u$ with $u\le \tilde{u}$ and  the proof of Theorem~\ref{thm:main} is complete.
$\hfill\Box$ 

\subsection{Proof of Theorem~\ref{thm:CS}}


\label{subsection:pf2}



Suppose $f$ satisfies {\rm (M)}, {\rm (S)}, {\rm (C)}, {\rm (L)} and is supercritical. Let $\varepsilon >0$ and $\gamma >0$ be as guaranteed by Theorem~\ref{thm:main}. 
For $\lam >0$, define the one-parameter family of nonnegative initial data $\mu^{\lam}\in\CB$ by
\[
  \mu^{\lam}(x)= 
  \psi\left( \lam |x|^{-\theta}\right),\qquad x\in\R^N.
\]
We first show that problem~\eqref{eq:SHE} is locally solvable with $ \mu=\mu^{\lam}$ for all $\lam$ small enough. 
By Theorem~\ref{thm:main}~(b), it is sufficient to show that there exists $r=r(\lam)>0$ such that $G(\mu^{\lam})\le \varepsilon |x|^{-\theta}$ for $|x|\le r$; i.e.,  
that   there  exists $y_0=y_0(\lam)>0$
such that
\begin{equation}
  \label{eq:v1}
G(\psi (y ))\le \frac{\varepsilon y}{\lam }\qquad\text{for all}\  y\ge y_0. 
\end{equation}

Now, since $\vp\asymp G$ and $G(u)\to\infty$ as $u\to\infty$, it follows that  $\vp (u)\to\infty$ as $u\to\infty$. Hence, since  $\vp\circ\psi\asymp\text{id}$,
we have  $\psi (u)\to\infty$ as $u\to\infty$. Thus, there exist $C_1, C_2>0$ and $y_0>0$ such that 
\begin{equation}
\label{eq:fdash1}
C_2\le \frac{G(\psi (y))}{y}= \frac{G(\psi (y))}{\vp(\psi(y))}
\frac{\vp(\psi(y))}{y}\le C_1
\qquad\text{for all}\  y\ge y_0.
\end{equation}
Setting $\lam_1=\varepsilon /C_1$, we see that
\eqref{eq:v1} holds for any $\lam\in (0,\lam_1)$ and therefore problem~\eqref{eq:SHE} is locally solvable with $\mu=\mu^{\lam}$ for all $\lam $ small enough.  

In a similar manner,  
we show that problem~\eqref{eq:SHE} is not locally solvable with $ \mu=\mu^{\lam}$ for large enough  $\lam$. 
By Theorem~\ref{thm:main}~(a), 
it is sufficient to show that $G(\mu^{\lam})\ge \gamma |x|^{-\theta}$  for  small enough $|x|$; i.e.,  that
\begin{equation}
  \label{eq:v2}
  G(\psi (y ))\ge \frac{\gamma y}{\lam }
\end{equation}
for large enough $y>0$. 
Setting $\lam_2=\gamma /C_2$, we see again from \eqref{eq:fdash1} that
\eqref{eq:v2} holds for all large enough $y >0$. Thus  problem~\eqref{eq:SHE} is not locally solvable with $ \mu=\mu^{\lam}$ for   large enough $\lam$.

Now set
\[
\lam_0=\sup\{\lam \ge 0:\ \text{problem~\eqref{eq:SHE} is locally solvable with  } \mu=\mu^{\lam}\}.
\]
By the arguments above, 
it follows that $\lam_0>0$ and finite. 
By definition of $\lam_0$,  problem~\eqref{eq:SHE} is  not locally solvable  for any $\lam >\lam_0$ and 
there exists a sequence $\lam_k\to\lam_0$ such that  problem~\eqref{eq:SHE} is locally solvable for $\lam =\lam_k$, for all $k\in\N$. 
Consequently, problem~\eqref{eq:SHE} is locally solvable for all $\lam \in (0,\lam_k)$, for all $k\in\N$. 
Thus problem~\eqref{eq:SHE} is locally solvable for all $\lam \in (0,\lam_0)$.  Taking $\mu_c=\mu^{\lam_0}$, we see that $\mu_c$ satisfies parts (i) and (ii) of Definition~\ref{def:CS}.

Finally we verify the local integrability of $\mu_c$. 
Since $f$ is supercritical, we may choose $\varepsilon_0>0$ small enough such that $N+\theta (1-q_f-\varepsilon_0)>0$. Then using \eqref{eq:4.5} in Lemma~\ref{Lemma:4.1}, and   $\psi(y)\le G^{-1}(C_1y)$ for large enough $y>0$ (by  \eqref{eq:fdash1}),   we  obtain
\begin{align*}
  \int_{B_r} \mu_c(x)\, \dee x 
  & = 
  \int_{B_r} \psi\left(  {\lam_0} |x|^{-\theta}\right) \dee x 
  \le   \int_{B_r} G^{-1}\left( {C}|x|^{-\theta}\right) \dee  x
  \\
  & =   
  C \int_0^{r} \sigma^{N-1}G^{-1}\left( {C}\sigma^{-\theta}\right) \dee \sigma
  = C \int_0^{r} \sigma^{N-1}F^{-1}\left(C\sigma^{\theta}\right) \dee \sigma 
  \\
  & \le   
  C 
  \int_0^{r} \sigma^{N-1+\theta (1-q_f-\varepsilon_0)}\,  \dee \sigma<\infty 
\end{align*}
for small enough $r>0$. 
Hence  $\mu_c\in  L^{1,+}_{{\rm loc}}\left({\mathbb R}^N\right)$, as required. 
It  follows that $\mu_c$ satisfies Definition~\ref{def:CS} and is a DCS of problem~\eqref{eq:SHE}.  The uniqueness of $\lam_0$ is an obvious consequence of its definition.  
$\hfill\Box$ 

\medskip




\section{Application: a generalised exponential}
\label{section:apps}

In this final section we apply our main results from Theorem~\ref{thm:main} and Theorem~\ref{thm:CS} to a family of 
 exponential-type nonlinearities which generalise problem~\eqref{eq:E} (in Section~\ref{subsec:background}) and Example~\ref{eg:exp}. 

We consider the problem
\begin{equation}
 \tag{GE}
 \label{eq:GE}
  \begin{cases}
    \partial_t u + (-\Delta)^\frac{\theta}{2}u =\exp_n (u^p), & x\in \mathbb{R}^N, \,\,\, t>0,
    \\[3pt]
    u(x,0) = \mu(x), & x\in \mathbb{R}^N,
  \end{cases}
\end{equation}
for any $n\in\N$ and $p>0$, 
where $\exp_n\colon [0,\infty )\to [\exp_n(0),\infty)$ denotes $n$ compositions of the exponential function; i.e.,   
\[
\exp_1(u) := e^u, 
\qquad 
\exp_{k+1}(u) := \exp\left(\exp_k(u)\right),\qquad k\in\N. 
\]
We denote by $\log_n\colon [\exp_n(0),\infty)\to [0,\infty )$ the inverse of the function $\exp_n$  
and write $e_n:=\exp_n(0)$. 
\begin{Theorem}
  \label{Theorem:GE}
  Let  $\mu\in L^{1,+}_{{\rm loc}}(\R^N)$, $n\in\N$ and $p>0$.
  Then the following hold:
  \begin{itemize}
    \item[{\rm (a)}]
    there exist $\gamma=\gamma(N,\theta,p)>0$ and $r>0$ such that if
    \begin{equation}\label{eq:5.2}
      \exp_n (\mu(x)^p) 
      \ge
     \frac{ \rule[-0.3cm]{0cm}{0cm} \gamma |x|^{-\theta}  \left( \log_n \left( |x|^{-1}\right) \right)^{\frac{1}{p}}}
      {\displaystyle\prod_{k=1}^{n} \log_k \left(  |x|^{-1} \right)}
    \end{equation}
    for a.a.~$x\in B_r$, 
    then problem~\eqref{eq:GE}  is not locally solvable;
    \item[{\rm (b)}] if $\mu\in\CB$ then 
    there exist $\varepsilon=\varepsilon(N,\theta,p)>0$ and $r>0$ such that if
    \begin{equation}\label{eq:5.3}
      \exp_n (\mu(x)^p) 
      \le
     \frac{ \rule[-0.3cm]{0cm}{0cm} \varepsilon |x|^{-\theta}  \left( \log_n \left( |x|^{-1}\right) \right)^{\frac{1}{p}}}
      {\displaystyle\prod_{k=1}^{n} \log_k \left(  |x|^{-1} \right)}
    \end{equation}
    for a.a.~$x\in B_r$, 
    then problem~\eqref{eq:GE}  is locally solvable;
    \item[{\rm (c)}] 
    there exist $r>0$ and $\lambda_0>0$ 
    such that problem~\eqref{eq:GE} possesses a dilation-critical singularity          
    given by 
    \begin{equation}
    \label{eq:newmu}
    {{\mu}_c}(x) =
    \left\{
    \begin{array}{ll}  
    \left(
     \log_n\left[ \frac{ \displaystyle{
    \rule[-0.5cm]{0cm}{0cm}
     \lambda_0 |x|^{-\theta}  \left( \log_n \left(\lambda_0 |x|^{-\theta}\right) \right)^{\frac{1}{p}}}}
    {\displaystyle{ \prod_{k=1}^{n} \log_k \left( \lambda_0 |x|^{-\theta} \right)}}
    \right]
    \right)^\frac{1}{p}, &  |x|<r,\\[35pt]
    0, & |x|\ge r.
    \end{array}
    \right.
    \end{equation}
  \end{itemize}
\end{Theorem}

\noindent{\textbf{Proof.}}
For $n\in \N$ and $p>0$, let  $f(u) = \exp_n(u^p)$. 
Then $f$ satisfies {\rm (M)}, {\rm (S)} and {\rm (C)}. 
Furthermore, thanks to l'H\^opital's rule (see Remark~\ref{rem:ORV}), we see that
\begin{align}
  \notag
 q_f= \lim_{u\to \infty} f'(u)F(u)
  =\lim_{u\to \infty} \frac{f'(u)^2}{f(u)f''(u)} 
  = 1. 
\end{align} 
{Hence $\pf =\infty$ and $f$ is supercritical.} Since $F_0<\infty$ and 
\begin{equation}
  \notag 
  f'(u) = pu^{p-1} \prod_{k=1}^{n} \exp_k(u^p),
\end{equation}
by Remark~\ref{rem:ORV}~(i) there exist $C_1,C_2>0$ such that  
\begin{equation}
  \label{eq:5.5}
  C_1(\max\{u,1\})^{p-1} \prod_{k=1}^{n} \exp_k(u^p)\le G(u) 
  \le C_2(\max\{u,1\})^{p-1} \prod_{k=1}^{n} \exp_k(u^p)
\end{equation}
for all $u\ge 0$.  

We first prove Theorem~\ref{Theorem:GE} part (a).  
Set 
\begin{equation} 
\label{eq:H}
H(u) = \prod_{k=1}^{n} \left( \log_k u \right)^{-1},\qquad  u>e_n. 
\end{equation} 
Let $\gamma\ge 1$, to be chosen later. 
Let $R\in (0,e_n^{-1})$ be such that  \eqref{eq:5.2} holds for all $x\in B_R$. 
For any $r\in (0,R]$, set 
\begin{equation}
  \label{eq:5.6}
  \hat{\mu}(x) = 
  \left[ \log_n \left(
    e_n + 
    \frac{\gamma}{2} |x|^{-\theta} \left( \log_n \left( e_n+ |x|^{-1} \right)  \right)^\frac{1}{p} 
    H\left( e_n+ |x|^{-1} \right)
    \chi_r(x)
  \right) \right]^\frac{1}{p},\qquad x\in\mathbb{R}^N. 
\end{equation}
By the nonnegativity of $\mu$ and \eqref{eq:5.2}, 
taking small enough $r\in(0,R]$ if necessary, 
we have
\begin{equation}
  \label{eq:5.8}
  \begin{split}
  \exp_n (\hat{\mu}(x)^p) 
   & = 
  e_n + 
  \frac{\gamma}{2} |x|^{-\theta}\left( \log_n \left(  e_n+  |x|^{-1} \right)   \right)^{\frac{1}{p}}
  H\left(   e_n+  |x|^{-1} \right)
  \chi_r(x)\\
  & \le 
  \begin{cases}
  \gamma |x|^{-\theta} \left( \log_n \left(|x|^{-1} \right)   \right)^\frac{1}{p}
  H\left(|x|^{-1} \right), & x\in B_r, \\ 
  e_n, & x\not\in B_r, 
  \end{cases} 
  \\
  & \le \exp_n (\mu(x)^p) 
 \end{split}
\end{equation}
for a.a.~$x\in{\mathbb R}^N$. 
Hence $\mu(x)\ge \hat{\mu}(x)$ for a.a.~$x\in\mathbb{R}^N$.
Therefore, by Lemma~\ref{Lemma:2.2}
it suffices to prove that problem~\eqref{eq:GE} is not locally solvable
for  $\mu=\hat{\mu}$. 

By \eqref{eq:5.6},  
taking large enough $\gamma$ and small enough $r\in(0,R)$ if necessary, 
we have  
\begin{equation*}
\begin{split}
  \exp_n(\hat{\mu}(x)^p)  
  & \ge C'\gamma|x|^{-\theta}
  \left( \log_n \left(   e_n+  |x|^{-1} \right)   \right)^\frac{1}{p}  
  H\left(   e_n+  |x|^{-1} \right) 
   \chi_r(x), \qquad  n\ge 1,\\ 
   \hat{\mu}(x) 
  & 
  \ge 
  C' \left( \log_n \left(   e_n+  |x|^{-1} \right)   \right)^\frac{1}{p} \chi_r(x), \qquad  n\ge 1,
\end{split}
\end{equation*}
and
\begin{equation*}
  \exp_k(\hat{\mu}(x)^p) 
   \ge 
  C' \left( \log_{n-k} \left(   e_n+  |x|^{-1} \right)   \right) \chi_r(x)\quad
  \mbox{ for $k\in \{1,\ldots,n-1\}$ and $n\ge 2$},
\end{equation*}
for $x\in{\mathbb R}^N$, where $C'$ is a positive constant independent of $\gamma$. 
These together with \eqref{eq:5.5}--\eqref{eq:5.8} imply that
\[
  G(\hat{\mu}(x))\ge C(\max\{\hat{\mu}(x),1\})^{p-1} 
  \prod_{k=1}^n \exp_k (\hat{\mu}(x)^p) 
  \ge 
  C \gamma |x|^{-\theta} \chi_r(x),\qquad x\in{\mathbb R}^N.
\]
Then, applying Theorem~\ref{thm:main}, 
we see that problem~\eqref{eq:GE} is not locally solvable if $\gamma$ is large enough. 
Thus Theorem~\ref{Theorem:GE}~(a) follows.

We now prove part (b) of Theorem~\ref{Theorem:GE}. 
Let $R\in (0,e_n^{-1})$ be such that \eqref{eq:5.3} holds for all $x\in B_R$. 
For any~$\varepsilon\in (0,1)$, set
\begin{equation}
  \label{eq:5.9}
  \tilde{\mu}(x) = 
  \left[ \log_n \left(
    K +  
    2\varepsilon |x|^{-\theta} \left( \log_n \left(  e_n+  |x|^{-1} \right)  \right)^\frac{1}{p}
    H\left(  e_n+  |x|^{-1} \right)  \chi_{r}(x) 
  \right) \right]^\frac{1}{p},\qquad x\in\mathbb{R}^N,
\end{equation}
where $K \ge e_n$ and $r\in (0,R]$.
Then, 
taking large enough $K$ and small enough $r$, 
by \eqref{eq:5.3} we have 
\[
\exp_n(\mu(x)^p)\le \exp_n(\tilde{\mu}(x)^p)
\] 
for a.a.~$x\in{\mathbb R}^N$, 
so that $\mu(x)\le \tilde{\mu}(x)$ for a.a.~$x\in\mathbb{R}^N$. 
Thanks to Lemma~\ref{Lemma:2.2}, 
it suffices to prove the existence of a solution of problem~\eqref{eq:GE} with~$\mu=\tilde{\mu}$.

By~\eqref{eq:5.9} we have 
\begin{equation}
  \notag 
  \begin{aligned}
    \exp_{n}(\tilde{\mu}(x)^p) 
    & = 
    K +  
    2\varepsilon |x|^{-\theta} \left( \log_n \left(  e_n+  |x|^{-1} \right)   \right)^\frac{1}{p}
    H\left(  e_n+  |x|^{-1} \right)  \chi_{r}(x) , \qquad  n\ge 1,
    \\
    \tilde{\mu}(x) 
    & \le 
    C \left(
      \log_{n} \left(  e_n+  |x|^{-1} \right)  
    \right)^\frac{1}{p} 
    \chi_{r}(x)  + C,  \qquad  n\ge 1,
  \end{aligned}
\end{equation}
and
\begin{equation*}
    \exp_{k}(\tilde{\mu}(x)^p) 
    \le 
    C \left(
      \log_{n-k}\left(  e_n+  |x|^{-1} \right) 
    \right)
    \chi_{r}(x) + C \quad 
\mbox{ for $k\in \{1,\ldots,n-1\}$ and $n\ge 2$}
\end{equation*}
for $x\in\mathbb{R}^N$.  
Then, 
by \eqref{eq:5.5} and \eqref{eq:H}
we obtain 
\begin{equation}
  \notag 
  G(\tilde{\mu}(x)) 
  \le 
  C(\max\{\tilde{\mu}(x),1\})^{p-1}
  \prod_{k=1}^n \exp_k(\tilde{\mu}(x)^p)
  \le 
  C \varepsilon |x|^{-\theta} + C 
\end{equation}
for $x\in\mathbb{R}^N$. 
Therefore, taking a small enough $\varepsilon>0$, we may apply Theorem~\ref{thm:main} 
to obtain a nonnegative solution of problem~\eqref{eq:GE} with $\mu=\tilde{\mu}$, as required.
To prove Theorem~\ref{Theorem:GE}~(c), we use Theorem~\ref{thm:CS}.
For $u>L$, set 
\begin{align}
    \label{eq:psi}
    A(u) = u\left(\log_n u\right)^\frac{1}{p}H(u), \quad\,\,\, 
    \varphi_0 (u) = u^{p-1}\prod_{k=1}^n \exp_k(u^p), \quad\,\,\, 
    \psi_0 (u)  = \left( \log_n
    A(u)
    \right)^\frac{1}{p},
\end{align}
where  $L>1$ is  chosen large enough so that $\vp_0$ and $\psi_0$ are both positive, increasing functions. 
Clearly $\vp_0\asymp f'$, so if we can show that $\vp_0\circ\psi_0 \asymp \text{id}$ 
then we may extend each of $\vp_0$ and $\psi_0$ by zero to have domain $[0,\infty )$ and use \eqref{eq:critsing}.

By \eqref{eq:psi} we have
\begin{equation}
\label{eq:comp}
\begin{split}
\vp_0(\psi_0(u)) & = 
\left(\log_n A (u)\right)^{\frac{p-1}{p}}\prod_{k=1}^{n} \exp_k \left(  \log_n A (u)\right)
=  \left( \log_n A (u)\right)^\frac{p-1}{p}\prod_{i=0}^{n-1} \log_i  A (u) \\
& = \left( \log_n A (u)\right)^{-\frac{1}{p}} A (u)\prod_{i=1}^{n}\log_i A (u)
=  u \left( \frac{\log_n A (u)}{\log_n u}\right)^{-\frac{1}{p}}\prod_{i=1}^{n}\frac{\log_i A (u)}{\log_i u}, 
\end{split}
\end{equation}
where $\log_0 = \mathrm{id}$. 
Next we claim that, for all $i\ge 1$, 
\begin{equation}
\lim_{u\to\infty}\frac{\log_i A (u)}{\log_i u} =1.\label{eq:limlog}
\end{equation}
For $i=1$ we have
\begin{equation*}
\lim_{u\to\infty}\frac{\log A (u)}{\log u}  =\lim_{u\to\infty} \frac{\log u +\frac{1}{p}\log_{n+1}  u-\sum_{j=1}^n\log_{j+1}  u }{\log u}=1.
\end{equation*}
If \eqref{eq:limlog} is true for some $i\ge 1$, then
\begin{equation*}
\lim_{u\to\infty}\frac{\log_{i+1} A (u)}{\log_{i+1} u}  =\lim_{u\to\infty}
\frac{\log \left[(\log_{i} u)\left(\frac{\log_{i} A (u)}{\log_{i} u}\right)\right]}{\log\left( \log_{i}u\right)}  
=1
\end{equation*}
and  \eqref{eq:limlog} follows by induction. 
It follows  from  \eqref{eq:comp}  and \eqref{eq:limlog} that
\bs
\lim_{u\to\infty} \frac{\vp_0(\psi_0(u))}{u}=1
\es
and so $\vp_0\circ\psi_0 \asymp \text{id}$. 
Finally, we extend each of $\vp_0$ and $\psi_0$ by zero on $[0,L]$ to yield functions  $\vp ,\psi\colon [0,\infty )\to [0,\infty )$ satisfying the hypotheses of Theorem~\ref{thm:CS}. From \eqref{eq:critsing}, \eqref{eq:H} and \eqref{eq:psi}  we obtain the DCS
\[
{\mu}_c(x)=\psi \left(\lam_0 |x|^{-\theta}\right),
\]
yielding \eqref{eq:newmu} with $r=(\lam_0/L)^{1/\theta}$. 
$\hfill\Box$

\begin{Remark}
  \label{Remark:4.1} 
  {\rm 
  Consider problem~\eqref{eq:GE} in the special case  $n=1$.
  For any $\lam >0$, let $\mu_\lam$ be the function in ${\mathbb R}^N$ satisfying 
  \[
  \exp(\mu_\lam^p)=\lam|x|^{-\theta} |\log |x||^{-\frac{(p-1)}{p}}
  \chi_{B_{r}}(x) + \chi_{\mathbb{R}^N\backslash  B_{r}}(x),\qquad x\in{\mathbb R}^N. 
  \]
Since 
  \[
  \mu_\lam(x)^{p-1}\exp(\mu_\lam^p)=\lam
  \theta^\frac{(p-1)}{p}
  |x|^{-\theta}(1+o(1))\qquad\mbox{as}\qquad x\to 0,
  \]
 we see that 
  \[
  (\max\{\mu_\lam,1\})^{p-1}\exp(\mu_\lam^p)\in L_{\rm uloc}^r\left({\mathbb R}^N\right)\qquad
  \mbox{if and only if}\qquad 1\le r<\frac{N}{\theta}.
  \]
Hence, in the framework of \cites{FI, GM}, it  is not possible to determine the existence of solutions with $\mu=\mu_\lambda$, 
even in the case of $p\ge 1$. Likewise, the arguments in \cites{FI, GM} are insufficient to establish the nonexistence of solutions for $\mu=\mu_\lambda$ and large $\lambda$. 
(Recall assertions~(E1), (E2) and (E3) in Section~\ref{subsec:background}.)

However, the theory we have developed here {\it does} allow one to determine existence and nonexistence results for this family of initial data. Specifically, by Theorem~\ref{Theorem:GE} (a)--(b), 
  we see that: 
  \begin{itemize}
    \item 
   problem~\eqref{eq:GE} with $n=1$ is not locally solvable   
    with $\mu=\mu_\lam$ if $\lam$ is large enough;
    \item 
   problem~\eqref{eq:GE} with $n=1$ is locally solvable 
    with $\mu=\mu_\lam$ if $\lam$ is small enough.
  \end{itemize}
%
  }
\end{Remark}

\medskip

\noindent{\bf Acknowledgments.}
YF was supported in part  by JSPS KAKENHI Grant Number 23K03179.
KH and KI  were supported in part by JSPS KAKENHI Grant Number JP19H05599.
All the authors  were supported in part by a Daiwa Anglo-Japanese Foundation Award (Ref: 14353/15194).
The authors would like to thank the reviewers for their careful attention and helpful suggestions.


\end{document}